\documentclass{article}
\usepackage{amsthm,amsmath,amssymb}
\usepackage{cite}
\usepackage{pstricks,pst-node}

\newtheorem{theorem}{Theorem}
\newtheorem{proposition}{Proposition}
\newtheorem{lemma}{Lemma}

\newcommand{\C}{\mathbb{C}}
\newcommand{\Z}{\mathbb{Z}}
\newcommand{\N}{\mathbb{N}}

\newcommand{\R}{\mathbb{R}}
\newcommand{\g}{\mathfrak{g}}
\newcommand{\heart}{\heartsuit}
\newcommand{\tensor}{\otimes}
\renewcommand{\bar}{\overline}
\renewcommand{\sp}{\operatorname{sp}}
\renewcommand{\sl}{\operatorname{sl}}
\DeclareMathOperator{\Inv}{Inv}
\DeclareMathOperator{\Hom}{Hom}
\DeclareMathOperator{\so}{so}
\newcommand{\ie}{{\em i.e.}}
\newcommand{\Ie}{{\em I.e.}}
\newcommand{\Thetabar}{\bar{\Theta}}

\psset{linewidth=.4pt,dash=3pt 3pt,doublesep=.05,dotsize=1pt 5}
\SpecialCoor

\newcommand{\middlearrow}{\lput{:U}{\pspicture(0,0)(0,0)
\psline[arrows=->,arrowscale=1.5](2.2pt,0)(2.3pt,0)\endpspicture}}
\newcommand{\redmiddlearrow}{\lput{:U}{\pspicture(0,0)(0,0)
\psline[arrows=->,arrowscale=1.5](2.2pt,0)(2.3pt,0)
\endpspicture}}
\newcommand{\unmiddlearrow}{\lput{:U}{\pspicture(0,0)(0,0)
\psline[arrows=->,arrowscale=1.5](-2.2pt,0)(-2.3pt,0)\endpspicture}}
\newcommand{\unredmiddlearrow}{\lput{:U}{\pspicture(0,0)(0,0)
\psline[arrows=->,arrowscale=1.5](-2.2pt,0)(-2.3pt,0)
\endpspicture}}

\newcommand{\littley}{%
\qline(.2887;30)(0,0)\qline(.2887;150)(0,0)\qline(.2887;270)(0,0)}
\newcommand{\littlelam}{%
\qline(.2887;90)(0,0)\qline(.2887;210)(0,0)\qline(.2887;330)(0,0)}

\newcommand{\psgoesto}{\hspace{.5cm}\pspicture[.5](0,-.1)(1,.1)
    \psline{->}(0,0)(1,0)
    \endpspicture\hspace{.5cm}}

\begin{document}
\title{Web bases for $\sl(3)$ are not dual canonical}
\author{Mikhail Khovanov
\thanks{Supported by a Sloan graduate fellowship in mathematics} \and
Greg Kuperberg \thanks{Partly supported by NSF grant \#9423300}}
\date{June 3, 1996}
\maketitle

\section{Introduction}

Given a simple Lie algebra $\g$ over $\C$ and a finite list of
finite-dimensional, irreducible representations $V_1,V_2,\ldots,V_n$,
one can study different bases of the tensor product representation
$$V_1 \tensor V_2 \tensor \ldots \tensor V_n$$
or its invariant space 
$$\Inv(V_1 \tensor V_2 \tensor \ldots \tensor V_n).$$
The quantum group $U_q(\g)$ has representations and vector spaces of
invariants which generalize these, and one can also study their bases,
with or without the intention of specializing to $q=1$.  (For
simplicity, we will usually consider $U_q(\g)$ as an algebra over
$\C(q^{1/2})$, and we will only occassionally mention $\Z[q^{\pm 1/2}]$
as a ground ring.)  Lusztig's remarkable canonical bases
\cite{Lusztig:book}, which are the same as Kashiwara's crystal bases
\cite{Kashiwara:crystallizing}, extend to bases of these spaces and have many
important properties.

When $\g = \sl(2)$, the Temperley-Lieb category \cite{KL:recoupling,FK:canonical} gives
another set of bases for the invariant spaces.  It was recently
established that these bases are dual canonical, \ie, dual in the sense
of linear algebra to canonical bases\cite{FK:canonical}.  The Temperley-Lieb
category gives a particularly explicit, simple, and useful definition
of the dual canonical bases of invariants (dual canonical invariants)
which establishes further natural properties of these bases.

Reference~\citen{Kuperberg:spiders} defines generalizations of the
Temperley-Lieb category to the three rank two Lie algebras $A_2 \cong
\sl(3)$, $B_2 \cong \sp(4) \cong \so(5)$, and $G_2$.  These
generalizations are called combinatorial rank two spiders.  The bases
they yield are called web bases and their individual basis vectors are
called webs.  One may conjecture that these bases are also dual
canonical.  As evidence for the conjecture, consider the following
properties which the $A_2$ web bases share with dual canonical
invariants:

\begin{description}
\item[1.] Let $V_1,V_2,\ldots,V_n$ be arbitrary irreducible
representations of $U_q(\sl(3))$.  Then there is a natural cyclic
permutation operator
$$\Inv(V_1 \tensor V_2 \tensor \ldots \tensor V_n) \to 
\Inv(V_2 \tensor V_3 \tensor \ldots \tensor V_n \tensor V_1),$$
and it sends basis webs to other basis webs.
(Reference~\citen{Lusztig:book}, Prop. 28.2.4, establishes this property for
dual canonical invariants.)
\item[2.] The tensor product of two basis webs is a basis web.  (For
dual canonical invariants, this is a corollary of
Theorem~\ref{thlusztig}.)
\item[3.] If two adjacent tensor factors of a basis web are
dual 3-dimensional representations, then contracting them
produces a linear combination of basis webs with coefficients in
$\N[-q^{1/2},-q^{-1/2}]$.  (Here $\N$ means the non-negative integers;
the property is conjectural for dual canonical bases.)
\item[4.] Let $V$ and $V'$ be tensor products of arbitrary irreducible
representations, and let $V(\lambda)$ be the irreducible representation
of highest weight $\lambda$.  Then $\Inv(V \tensor B)$ decomposes as
$$\Inv(A \tensor B) \cong \bigoplus_\lambda \Inv(A \tensor V(\lambda))
\tensor \Inv(V(\lambda^*) \tensor B).$$
This decomposition induces a grading by $\lambda$, which leads to two
filtrations by the usual partial ordering on dominant weights.  The web
basis refines the ascending filtration.

\item[5.] The web bases are dual canonical in small cases.

\end{description}

In this paper, we will disprove the conjecture. Let $V^+$ be the defining
3-dimensional representation of $\sl(3)$ and let $V^-$ be the dual
representation, and let $V^+$ and $V^-$ also denote the corresponding
representations of $U_q(\sl(3))$. Then:

\begin{theorem} Every basis web in
$$\Inv(V_1 \tensor V_2 \tensor \ldots \tensor V_n),$$
where each $V_i$ is either $V^+$ or $V^-$, is dual canonical when
$n \le 12$, except for a single basis web in
$$\Inv((V^+ \tensor V^+ \tensor V^- \tensor V^-)^{\tensor 3})$$
and its counterparts given by cyclic permutation of tensor factors. 
\label{thmain}
\end{theorem}

To see the extent of early agreement between the two kinds of bases,
note that there are 35 permutations of six tensor factors of $V^+$ and
six tensor factors of $V^-$ which are inequivalent under sign flip,
reversal of order, and cyclic permutation.  Each permutation yields a
513-dimensional vector space of invariants. All 513 basis webs are dual
canonical unless the $V^+$'s and the $V^-$'s are in the arrangement
stated in the theorem, in which case 512 of them are. However, the
fraction of basis webs that are dual canonical must go to 0
exponentially as $n \to \infty$.

In comparing the two types of bases, we will often refer to the book by
Lusztig \cite{Lusztig:book}.  The results cited there are stated in terms of
canonical bases, but they can be translated to statements about dual
canonical bases.

\subsection{Acknowledgements}

The authors would like to thank Igor Frenkel for his attention
to this work.

\section{The quantum group $U_q(\sl(3))$}

The quantum group $U_q(\sl(3))$ is an associative algebra over
$\C(q^{1/2})$, where $q^{1/2}$ is an indeterminate, with generators
$E_i$, $F_i$, $K_i$, and $K_i^{-1}$ for $i=1,2$, and the following
relations:
\begin{gather*}
K_i K_i^{-1} =  K_i^{-1}K_i = 1\\
K_i K_j = K_j K_i \\
K_i E_j = q^{a_{ij}/2} E_j K_i \\
K_i F_j = q^{-a_{ij}/2}F_j K_i \\
E_i F_i - F_i E_i = \delta_{ij} \frac{K_i-K_i^{-1}}{q-q^{-1}} \\
E_i^2 E_j - [2] E_i E_j E_i + E_j E_i^2 = 0 \quad i \ne j \\
F_i^2 F_j - [2] F_i F_j F_i + F_j F_i^2 = 0 \quad i \ne j
\end{gather*}
Here $\delta_{ij}$ is 1 when $i=j$ and 0 when $i \ne j$, while $a_{ij}
= 3\delta_{ij} - 1$ is the Cartan matrix of $\sl(3)$.  The quantity
$[n]$ is a quantum integer, defined by
$$[n] = \frac{q^{n/2} - q^{-n/2}}{q^{1/2} - q^{-1/2}}.$$

After clearing denominators in the relations, one obtains a Hopf
algebra over $\Z[q^{1/2},q^{-1/2}]$.  For convenience, let
$v = -q^{1/2}$.  (Our $v$ is the negative of the $v$ in
Reference~\citen{Lusztig:book}.)

The algebra $U_q(\sl(3))$ is a Hopf algebra with a certain standard
coproduct $\Delta$.  In this paper, we will use a second coproduct
$\bar{\Delta}$ which is more appropriate for dual canonical bases.
This coproduct takes the following values on generators:
\begin{align*}
\bar{\Delta}(K_i^{\pm 1}) &= K_i^{\pm 1} \tensor K_i^{\pm 1} \\
\bar{\Delta}(E_i) &= E_i \tensor 1 + K_i^{-1} \tensor E_i \\
\bar{\Delta}(F_i) &= F_i \tensor K_i + 1 \tensor F_i
\end{align*}
We use this coproduct to understand $V \tensor V'$ as a representation of
$U_q(\sl(3))$ if $V$ and $V'$ are themselves representations.  Also, in any
representation, we will say that $e$ is an invariant vector if $Xe =
\epsilon(X)e$, where $\epsilon$ is a homomorphism from $U_q(\sl(3))$ to
$\C(v)$ given on generators by
\begin{align*}
\epsilon(E_i) &= \epsilon(F_i) = 0 \\
\epsilon(K_i) &= 1
\end{align*}
The vector space of all invariants of $V$ is denoted $\Inv(V)$.

The two irreducible representations of the quantum group $U_q(\sl(3))$ that we
will study are the 3-dimensional representations $V^+$ and $V^-$. We choose a
basis $e^\pm_{-1},e^\pm_{0}, e^\pm_1$ of $V^\pm$; the action of $U_q(\sl(3))$
on $V^+$ is given by:
\begin{align*}
K_1(e^+_1) &= q^{1/2} e^+_1 & K_1(e^+_0) &= q^{-1/2}e^+_0
& K_1(e^+_{-1}) &= e^+_{-1} \\
K_2(e^+_1) &= e^+_1  & K_2(e^+_0) &= q^{1/2}e^+_0
& K_2(e^+_{-1}) &= q^{-1/2}e^+_{-1} \\
E_1(e^+_0) &= e^+_1    & F_1(e^+_1) &= e^+_0 \\
E_2(e^+_{-1}) &= e^+_0 & F_2(e^+_0) &= e^+_{-1}
\end{align*}
and all other combinations such as $E_1(v_1)$ are 0.  Similarly,
the action on $V^-$ is given by:
\begin{align*}
K_1(e^-_1) &= e^-_1  & K_1(e^-_0) &= q^{1/2}e^-_0
& K_1(e^-_{-1}) &= q^{-1/2}e^-_{-1} \\
K_2(e^-_1) &= q^{1/2}e^-_1 & K_2(e^-_0) &= q^{-1/2}e^-_0
& K_2(e^-_{-1}) &= e^-_{-1} \\
E_1(e^-_{-1}) &= e^-_0 & F_1(e^-_0) &= e^-_{-1} \\
E_2(e^-_0) &= e^-_1    & F_2(e^-_1) &= e^-_0
\end{align*}
and all other combinations are 0.  These actions are
summarized by the weight diagrams of $V^+$ and $V^-$:
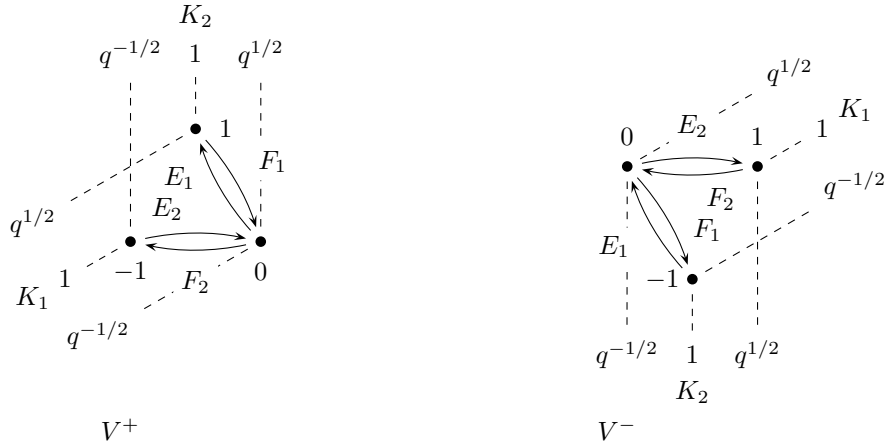
$$
\pspicture(-2.8,-3.3)(2.8,2.8)
\pnode(1;90){a1}\pnode(1;210){a2}\pnode(1;330){a3}
\qdisk(a1){2pt}\qdisk(a2){2pt}\qdisk(a3){2pt}
\rput(.4,1){$1$}\rput(-.866,-.9){$-1$}\rput(.866,-.9){$0$}
\rput(-1.732,-1){\rnode{k12}{$1$}}
\rput(-2.165,-.25){\rnode{k11}{$q^{1/2}$}}
\rput(-1.299,-1.75){\rnode{k13}{$q^{-1/2}$}}\rput(-2.165,-1.25){$K_1$}
\ncline[nodesep=.2,linestyle=dashed]{k11}{a1}
\ncline[nodesep=.2,linestyle=dashed]{k12}{a2}
\ncline[nodesep=.2,linestyle=dashed]{k13}{a3}
\rput(0,2){\rnode{k22}{$1$}}
\rput(-.866,2){\rnode{k21}{$q^{-1/2}$}}
\rput(.866,2){\rnode{k23}{$q^{1/2}$}}\rput(0,2.5){$K_2$}
\ncline[nodesep=.2,linestyle=dashed]{k21}{a2}
\ncline[nodesep=.2,linestyle=dashed]{k22}{a1}
\ncline[nodesep=.2,linestyle=dashed]{k23}{a3}
\ncarc[nodesep=.2,arrows=->,arcangle=12]{a2}{a3}\aput(.25){$E_2$}
\ncarc[nodesep=.2,arrows=->,arcangle=12]{a3}{a2}\Aput{\psframebox*{$F_2$}}
\ncarc[nodesep=.2,arrows=->,arcangle=12]{a3}{a1}\aput(.75){$E_1$}
\ncarc[nodesep=.2,arrows=->,arcangle=12]{a1}{a3}\Aput{\psframebox*{$F_1$}}
\rput(-1,-3){$V^+$}
\endpspicture\hspace{1cm}
\pspicture(-2.8,-3.3)(2.8,2.8)
\pnode(-1;90){a1}\pnode(-1;210){a2}\pnode(-1;330){a3}
\qdisk(a1){2pt}\qdisk(a2){2pt}\qdisk(a3){2pt}
\rput(-.4,-1){$-1$}\rput(.866,.9){$1$}\rput(-.866,.9){$0$}
\rput(1.732,1){\rnode{k12}{$1$}}
\rput(2.165,.25){\rnode{k11}{$q^{-1/2}$}}
\rput(1.299,1.75){\rnode{k13}{$q^{1/2}$}}\rput(2.165,1.25){$K_1$}
\ncline[nodesep=.2,linestyle=dashed]{k11}{a1}
\ncline[nodesep=.2,linestyle=dashed]{k12}{a2}
\ncline[nodesep=.2,linestyle=dashed]{k13}{a3}
\rput(0,-2){\rnode{k22}{$1$}}
\rput(.866,-2){\rnode{k21}{$q^{1/2}$}}
\rput(-.866,-2){\rnode{k23}{$q^{-1/2}$}}\rput(0,-2.5){$K_2$}
\ncline[nodesep=.2,linestyle=dashed]{k21}{a2}
\ncline[nodesep=.2,linestyle=dashed]{k22}{a1}
\ncline[nodesep=.2,linestyle=dashed]{k23}{a3}
\ncarc[nodesep=.2,arrows=->,arcangle=12]{a2}{a3}\aput(.25){$F_2$}
\ncarc[nodesep=.2,arrows=->,arcangle=12]{a3}{a2}\Aput{\psframebox*{$E_2$}}
\ncarc[nodesep=.2,arrows=->,arcangle=12]{a3}{a1}\aput(.75){$F_1$}
\ncarc[nodesep=.2,arrows=->,arcangle=12]{a1}{a3}\Aput{\psframebox*{$E_1$}}
\rput(-1,-3){$V^-$}
\endpspicture
$$
These bases for $V^+$ and $V^-$ are dual canonical because
they satisfy the conditions of bases at $\infty$ \cite[20.1.1]{Lusztig:book}.

\section{The combinatorial $A_2$ spider}

Strictly speaking, the combinatorial $A_2$ spider is an abstract tensor
category given by generators and relations which is known to be
isomorphic to the category of intertwiners of $U_q(\sl(3))$
\cite{Kuperberg:spiders}.  However, in this paper, this isomorphism
will be implicit and we will instead describe vectors
of $U_q(\sl(3))$, called webs, which are associated to it.

Following Reshetikhin and Turaev \cite{RT:ribbon}, we denote invariants and
equivariants in the representation category of $U_q(\sl(3))$ by means of
planar graphs.  If $V$ is a tensor product of $V^\pm$'s, an element of
$\Inv(V)$ may be denoted by some graph with oriented edges, with vertices
labelled by invariant tensors, and with a univalent vertex for each tensor
factor.  For each factor of $V^+$, the edge is oriented towards the
corresponding vertex, and for each factor of $V^-$, the edge is oriented
away:
$$
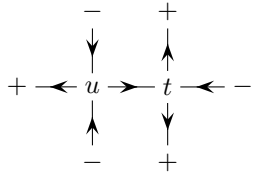
\pspicture(-2,-1.1)(2,1.1)
\rput( .5, 0){\rnode{a1}{$t$}}\rput(-.5, 0){\rnode{a2}{$u$}}
\rput( .5, 1){\rnode{b1}{$+$}}\rput( .5,-1){\rnode{b2}{$+$}}
\rput(1.5, 0){\rnode{b3}{$-$}}\rput(-.5, 1){\rnode{b4}{$-$}}
\rput(-.5,-1){\rnode{b5}{$-$}}\rput(-1.5, 0){\rnode{b6}{$+$}}
\ncline[nodesep=3pt]{a1}{b1}\middlearrow
\ncline[nodesep=3pt]{a1}{b2}\middlearrow
\ncline[nodesep=3pt]{b3}{a1}\middlearrow
\ncline[nodesep=3pt]{b4}{a2}\middlearrow
\ncline[nodesep=3pt]{b5}{a2}\middlearrow
\ncline[nodesep=3pt]{a2}{b6}\middlearrow
\ncline[nodesep=3pt]{a2}{a1}\middlearrow
\endpspicture
$$
As in this example, we abbreviate $V^+$ and $V^-$ by their signs.
Sometimes the signs or the orientations or both will be omitted in
cases where they are irrelevant or clear from context.

If $V$ and $V'$ are two different tensor products, 
the equivalence $\Hom(V,V') \cong \Inv(V^* \tensor V')$ will also be
important.  For example, a graph such as:
$$
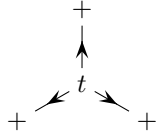
\pspicture(-1.2,-1.2)(1.2,1.2)
\rput(1;90){\rnode{a1}{$+$}}\rput(1;210){\rnode{a2}{$+$}}
\rput(1;330){\rnode{a3}{$+$}}\rput(0,0){\rnode{b}{$t$}}
\ncline[nodesep=3pt]{b}{a1}\middlearrow
\ncline[nodesep=3pt]{b}{a2}\middlearrow
\ncline[nodesep=3pt]{b}{a3}\middlearrow
\endpspicture
$$
might denote an element of $\Hom(V^-,V^+ \tensor V^+)$ just as well as
an element of $\Inv((V^+)^{\tensor 3})$. Using this equivalence,
compositions of homomorphisms are in general denoted by concatenation
and tensor products of homomorphisms (or invariants) are denoted by
juxtaposition, or disjoint union.

The combinatorial $A_2$ webs can be constructed from four elementary
invariants and the operations of tensor product and contraction.
The four invariants are
\begin{align*}
b^{+-} &= e^+_1 \tensor e^-_{-1} + v^{-1} e^+_0 \tensor e^-_0
+ v^{-2} e^+_{-1} \tensor e^-_{1} \\
b^{-+} &= e^-_1 \tensor e^+_{-1} + v^{-1} e^-_0 \tensor e^+_0
+ v^{-2} e^-_{-1} \tensor e^+_{1} \\
t^{---} &= e^-_1 \tensor e^-_0 \tensor e^-_{-1}
+ v^{-1} e^-_0 \tensor e^-_1 \tensor e^-_{-1}
+ v^{-1} e^-_1 \tensor e^-_{-1} \tensor e^-_0 \\
&\qquad + v^{-2} e^-_0 \tensor e^-_{-1} \tensor e^-_1
+ v^{-2} e^-_{-1} \tensor e^-_1 \tensor e^-_0
+ v^{-3} e^-_{-1} \tensor e^-_0 \tensor e^-_1 \\
t^{+++} &= e^+_1 \tensor e^+_0 \tensor e^+_{-1}
+ v^{-1} e^+_0 \tensor e^+_1 \tensor e^+_{-1}
+ v^{-1} e^+_1 \tensor e^+_{-1} \tensor e^+_0 \\
&\qquad + v^{-2} e^+_0 \tensor e^+_{-1} \tensor e^+_1
+ v^{-2} e^+_{-1} \tensor e^+_1 \tensor e^+_0
+ v^{-3} e^+_{-1} \tensor e^+_0 \tensor e^+_1
\end{align*}
The contraction operations are defined using the equivariant pairings
$\sigma_{+-}:V^+ \tensor V^- \to \C(q)$ and $\sigma_{-+}:V^- \tensor
V^+ \to \C(q)$, which are given by
\begin{alignat*}{2}
\sigma_{+-}(e^+_{-1} \tensor e^-_1) &=&
\sigma_{-+}(e^-_{-1} \tensor e^+_1) &= 1 \\
\sigma_{+-}(e^+_0 \tensor e^-_0)    &=&
\sigma_{-+}(e^-_0 \tensor e^+_0 )   &= v \\
\sigma_{+-}(e^+_1 \tensor e^-_{-1}) &=&
\sigma_{-+}(e^-_1 \tensor e^+_{-1}) &= v^2
\end{alignat*}
and all other values on basis vectors are 0.  Since the $\sigma$'s
are equivariant, they induce contractions
$$\Inv(V \tensor V^\pm \tensor V^\mp \tensor V') \to \Inv(V \tensor V')$$
for arbitrary tensor products $V$ and $V'$.

The planar graphs corresponding to the $b$'s and $t$'s are:
$$
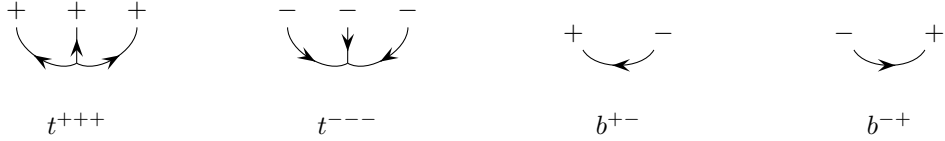
\pspicture(-.8,-.8)(.8,1.1)
\rput(-.8,1){\rnode{a1}{$+$}}\rput(0,1){\rnode{a2}{$+$}}
\rput(.8,1){\rnode{a3}{$+$}}\pnode(0,.3){b}
\nccurve[nodesepB=3pt,angleB=270,angleA=210]{b}{a1}\middlearrow
\nccurve[nodesepB=3pt,angleB=270,angleA=90]{b}{a2}\middlearrow
\nccurve[nodesepB=3pt,angleB=270,angleA=330]{b}{a3}\middlearrow
\rput(0,-.5){$t^{+++}$}
\endpspicture
\hspace{2cm}
\pspicture(-.8,-.8)(.8,1.1)
\rput(-.8,1){\rnode{a1}{$-$}}\rput(0,1){\rnode{a2}{$-$}}
\rput(.8,1){\rnode{a3}{$-$}}\pnode(0,.3){b}
\nccurve[nodesepA=3pt,angleA=270,angleB=210]{a1}{b}\middlearrow
\nccurve[nodesepA=3pt,angleA=270,angleB=90]{a2}{b}\middlearrow
\nccurve[nodesepA=3pt,angleA=270,angleB=330]{a3}{b}\middlearrow
\rput(0,-.5){$t^{---}$}
\endpspicture
\hspace{2cm}
\pspicture(-.8,-.8)(.8,1.1)
\rput(-.6,.7){\rnode{b1}{$+$}}\rput(.6,.7){\rnode{b2}{$-$}}
\ncarc[arcangle=60,nodesep=3pt]{b2}{b1}\middlearrow
\rput(0,-.5){$b^{+-}$}
\endpspicture
\hspace{2cm}
\pspicture(-.8,-.8)(.8,1.1)
\rput(-.6,.7){\rnode{b1}{$-$}}\rput(.6,.7){\rnode{b2}{$+$}}
\ncarc[arcangle=-60,nodesep=3pt]{b1}{b2}\middlearrow
\rput(0,-.5){$b^{-+}$}
\endpspicture
$$
while those corresponding to the $\sigma$'s are
$$
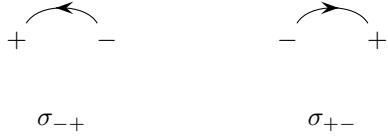
\pspicture(-.8,-.8)(.8,1.1)
\rput(-.6,.6){\rnode{b1}{$+$}}\rput(.6,.6){\rnode{b2}{$-$}}
\ncarc[arcangle=-60,nodesep=3pt]{b2}{b1}\middlearrow
\rput(0,-.5){$\sigma_{-+}$}
\endpspicture
\hspace{2cm}
\pspicture(-.8,-.8)(.8,1.1)
\rput(-.6,.6){\rnode{b1}{$-$}}\rput(.6,.6){\rnode{b2}{$+$}}
\ncarc[arcangle=60,nodesep=3pt]{b1}{b2}\middlearrow
\rput(0,-.5){$\sigma_{+-}$}
\endpspicture
$$

If the $b$'s and $t$'s are understood as equivariant homomorphisms
from the ground field $\C(v)$ to the corresponding invariant spaces,
they and the $\sigma$'s can be composed to form planar graphs.  For
example, the graph
$$
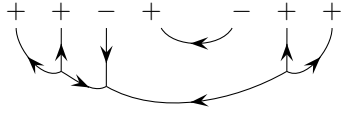
\pspicture(-.2,-2)(4.4,.2)
\rput(0,0){\rnode{a1}{$+$}}\rput(.6,0){\rnode{a2}{$+$}}
\rput(1.2,0){\rnode{a3}{$-$}}\rput(1.8,0){\rnode{a4}{$+$}}
\rput(3.0,0){\rnode{a5}{$-$}}\rput(3.6,0){\rnode{a6}{$+$}}
\rput(4.2,0){\rnode{a7}{$+$}}
\pnode(.6,-.8){b1}\pnode(1.2,-1){b2}\pnode(3.6,-.8){b3}
\nccurve[nodesepA=3pt,angleA=270,angleB=210]{a1}{b1}\unmiddlearrow
\nccurve[nodesepA=3pt,angleA=270,angleB=90]{a2}{b1}\unmiddlearrow
\nccurve[nodesepA=3pt,angleA=270,angleB=90]{a3}{b2}\middlearrow
\nccurve[nodesep=3pt,angleA=-60,angleB=240]{a4}{a5}\unmiddlearrow
\nccurve[nodesepA=3pt,angleA=270,angleB=90]{a6}{b3}\unmiddlearrow
\nccurve[nodesepA=3pt,angleA=270,angleB=330]{a7}{b3}\unmiddlearrow
\nccurve[angleA=330,angleB=210]{b1}{b2}\middlearrow
\nccurve[angleA=330,angleB=210]{b2}{b3}\unmiddlearrow
\endpspicture
$$
denotes the tensor
$$(I \tensor I \tensor \sigma_{+-} \tensor I \tensor b^{+-} \tensor
\sigma_{-+} \tensor I \tensor I)\circ (t^{+++} \tensor t^{---}
\tensor t^{+++}).$$

Given the identities
\begin{gather}
(I \tensor \sigma_{+-})\circ (b^{-+} \tensor I)
= (\sigma_{-+} \tensor I)\circ (I \tensor b^{+-}) = I \nonumber \\
(I \tensor \sigma_{-+})\circ (b^{+-} \tensor I)
= (\sigma_{+-} \tensor I)\circ (I \tensor b^{-+}) = I, \label{eisotopy}
\end{gather}
the value of a planar graph as a tensor is invariant under
isotopy of the graph.

In the combinatorial $A_2$ spider, a (monomial) web is defined as any
composition of tensor products of $b$'s, $t$'s, and $\sigma$'s.  Any
such web is denoted by an oriented graph in a disk with trivalent and
univalent vertices, and possibly closed loops, such that the edges are
either all out or all in at the trivalent vertices, and such that the
univalent vertices are at the boundary of the disk.

By the fundamental theorem of invariant theory, the set of all monomial
webs with any given boundary spans the corresponding set of invariants.
Moreover, the following relations hold:
\begin{eqnarray}
\pspicture[.4](-.6,-.5)(.6,.5)
\pscircle(0,0){.4}\psline[arrows=->,arrowscale=1.5](.1,.4)(.11,.4)
\endpspicture
& = & [3] \nonumber \\
\pspicture[.45](-1.5,-.8)(1.5,.8)
\rput(-1.2,0){\rnode{a1}{$-$}}
\pnode(-.4,0){a2}\pnode(.4,0){a3}
\rput(1.2,0){\rnode{a4}{$+$}}
\ncline[nodesepA=3pt]{a1}{a2}\middlearrow
\nccurve[angleA=120,angleB=60,nodesep=.3pt]{a3}{a2}\middlearrow
\nccurve[angleA=-120,angleB=-60,nodesep=.3pt]{a3}{a2}\middlearrow
\ncline[nodesepB=3pt]{a3}{a4}\middlearrow
\endpspicture
& = & -[2] \pspicture[.45](-1.2,-.8)(.8,.8)
\rput(-.6,0){\rnode{a1}{$-$}}\rput(.6,0){\rnode{a2}{$+$}}
\ncline[nodesep=3pt]{a1}{a2}\middlearrow
\endpspicture \nonumber \\
\pspicture[.42](-1.5,-1.4)(1.5,1.4)
\pnode(.4; 45){a1}\rput(1.2; 45){\rnode{b1}{$-$}}
\pnode(.4;135){a2}\rput(1.2;135){\rnode{b2}{$+$}}
\pnode(.4;225){a3}\rput(1.2;225){\rnode{b3}{$-$}}
\pnode(.4;315){a4}\rput(1.2;315){\rnode{b4}{$+$}}
\ncline[nodesepA=3pt]{b1}{a1}\middlearrow
\ncline[nodesepB=3pt]{a2}{b2}\middlearrow
\ncline[nodesepA=3pt]{b3}{a3}\middlearrow
\ncline[nodesepB=3pt]{a4}{b4}\middlearrow
\ncarc[arcangle= 15]{a2}{a1}\middlearrow
\ncarc[arcangle=-15]{a2}{a3}\middlearrow
\ncarc[arcangle= 15]{a4}{a3}\middlearrow
\ncarc[arcangle=-15]{a4}{a1}\middlearrow
\endpspicture
& = &
\pspicture[.4](-1,-.9)(1,.9)
\rput(.7; 45){\rnode{a1}{$-$}}\rput(.7;135){\rnode{a2}{$+$}}
\rput(.7;225){\rnode{a3}{$-$}}\rput(.7;315){\rnode{a4}{$+$}}
\ncarc[arcangle=-45]{a3}{a2}\middlearrow
\ncarc[arcangle=-45]{a1}{a4}\middlearrow
\endpspicture
+
\pspicture[.4](-1,-.9)(1,.9)
\rput(.7; 45){\rnode{a1}{$+$}}\rput(.7;135){\rnode{a2}{$-$}}
\rput(.7;225){\rnode{a3}{$+$}}\rput(.7;315){\rnode{a4}{$-$}}
\ncarc[arcangle=-45]{a3}{a2}\middlearrow
\ncarc[arcangle=-45]{a1}{a4}\middlearrow
\endpspicture
\label{fa2elliptic}
\end{eqnarray}
Thus, the set of non-elliptic webs, \ie, webs such that all internal
faces have at least six sides, also spans.  It is a fundamental result
that the set of non-elliptic webs is a basis of each invariant space
\cite{Kuperberg:spiders}.  These are the web bases that we will compare
to the dual canonical bases.

\section{State sums}

In the given bases of $V^+$ and $V^-$, the $b$, $t$, and $\sigma$
tensors have matrices, and any monomial web can be evaluated by the
usual linear algebra method of summing over indices of these matrices.
Such an expansion is equivalent to a state sum in the sense of
statistical mechanics.  Given a monomial web $w$, we first divide each
edge into segments whose edges are the points where the edge has a
horizontal tangent.  A state is then a function from the segments to
the set $\{-1,0,1\}$.  The weight of a state at each trivalent vertex
or horizontal tangent is a matrix entry of the corresponding $b$, $t$,
or $\sigma$ tensor. A boundary state is a function from just those
segments with univalent vertices to $\{-1,0,1\}$.  The weight of a
boundary state is then defined as the total weight of all extensions of
the boundary state to a state of the entire graph. The weights of the
boundary state are then the coefficients of $w$ expanded in the
tensor product basis.

We will abbreviate a state in a state sum by flow lines.  A collection
of flow lines in a monomial web $w$ is a subgraph that contains exactly
two of the three edges incident to each trivalent vertex.  Each flow
line is oriented; this orientation need not agree with the orientation
of $w$.  Every segment of an edge disjoint from a flow line has the
state 0.  If a flow line is oriented downward along a segment, the
segment has the state 1, while if it is oriented upward, the segment
has the state -1.  In this way, flow lines represent precisely those
states with non-zero weight.

A second convenience for computing state sums is to introduce
the linear endomorphisms
\begin{align*}
t^{++}_-&:V^- \to V^+ \tensor V^+ \\
t^{--}_+&:V^+ \to V^- \tensor V^- \\
t^+_{--}&:V^- \tensor V^- \to V^+ \\
t^-_{++}&:V^+ \tensor V^+ \to V^-
\end{align*}
They may be defined as the webs:
$$
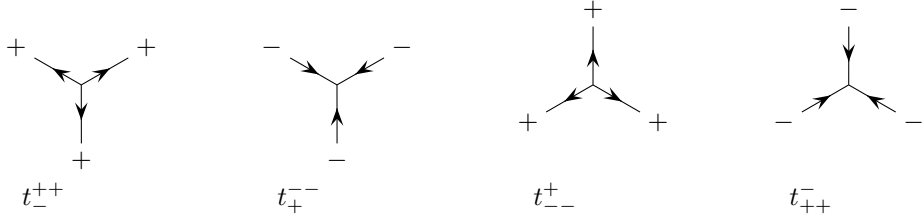
\pspicture(-1.2,-1.8)(1.2,1.2)
\rput(-1;90){\rnode{a1}{$+$}}\rput(-1;210){\rnode{a2}{$+$}}
\rput(-1;330){\rnode{a3}{$+$}}\pnode(0,0){b}
\ncline[nodesepB=3pt]{b}{a1}\middlearrow
\ncline[nodesepB=3pt]{b}{a2}\middlearrow
\ncline[nodesepB=3pt]{b}{a3}\middlearrow
\rput(-.5,-1.5){$t^{++}_-$}
\endpspicture\hspace{1cm}
\pspicture(-1.2,-1.8)(1.2,1.2)
\rput(-1;90){\rnode{a1}{$-$}}\rput(-1;210){\rnode{a2}{$-$}}
\rput(-1;330){\rnode{a3}{$-$}}\pnode(0,0){b}
\ncline[nodesepA=3pt]{a1}{b}\middlearrow
\ncline[nodesepA=3pt]{a2}{b}\middlearrow
\ncline[nodesepA=3pt]{a3}{b}\middlearrow
\rput(-.5,-1.5){$t^{--}_+$}
\endpspicture\hspace{1cm}
\pspicture(-1.2,-1.8)(1.2,1.2)
\rput(1;90){\rnode{a1}{$+$}}\rput(1;210){\rnode{a2}{$+$}}
\rput(1;330){\rnode{a3}{$+$}}\pnode(0,0){b}
\ncline[nodesepB=3pt]{b}{a1}\middlearrow
\ncline[nodesepB=3pt]{b}{a2}\middlearrow
\ncline[nodesepB=3pt]{b}{a3}\middlearrow
\rput(-.5,-1.5){$t^+_{--}$}
\endpspicture\hspace{1cm}
\pspicture(-1.2,-1.8)(1.2,1.2)
\rput(1;90){\rnode{a1}{$-$}}\rput(1;210){\rnode{a2}{$-$}}
\rput(1;330){\rnode{a3}{$-$}}\pnode(0,0){b}
\ncline[nodesepA=3pt]{a1}{b}\middlearrow
\ncline[nodesepA=3pt]{a2}{b}\middlearrow
\ncline[nodesepA=3pt]{a3}{b}\middlearrow
\rput(-.5,-1.5){$t^-_{++}$}
\endpspicture
$$
(Note that here, as before, the signs of the tensor subscripts are
opposite to the signs at the bottom of the webs, because
$\Hom(V,V') \cong \Inv(V^* \tensor V')$.)
Their coefficients are given by:
\begin{align*}
t^{++}_- e^-_1 &= e^+_1 \tensor e^+_0 + v^{-1} e^+_0 \tensor e^+_1 \\
t^{++}_- e^-_0 &= e^+_1 \tensor e^+_{-1} + v^{-1} e^+_{-1} \tensor e^+_1 \\
t^{++}_- e^-_{-1} &= e^+_0 \tensor e^+_{-1} + v^{-1} e^+_{-1} \tensor e^+_0
\end{align*}
and
\begin{align*}
t^+_{--} (e^-_1 \tensor e^-_0) &= v e^+_1 &
t^+_{--} (e^-_0 \tensor e^-_1) &= e^+_1 \\
t^+_{--} (e^-_1 \tensor e^-_{-1}) &= v e^+_0 &
t^+_{--} (e^-_{-1} \tensor e^-_1) &= e^+_0 \\
t^+_{--} (e^-_0 \tensor e^-_{-1}) &= v e^+_{-1} &
t^+_{--} (e^-_{-1} \tensor e^-_0) &= e^+_{-1}
\end{align*}
As before, all combinations not listed are 0.  The formulas for
$t^-_{++}$ and $t^{--}_+$ are the same; one just switches $+$'s and
$-$'s.  These endomorphisms have graphs that are Y's and $\lambda$'s;
their weights may be abbreviated with flow lines according to the
following chart:
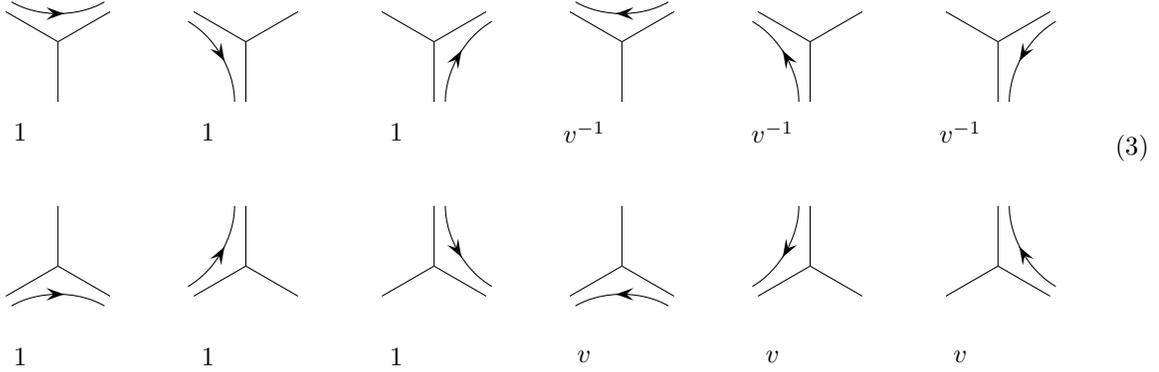
\begin{equation}
\pspicture(-1,-1.5)(1,1)
\pnode(-.8;90){a1}\pnode(-.8;210){a2}\pnode(-.8;330){a3}\pnode(0,0){b}
\ncline{a1}{b}\ncline{a2}{b}\ncline{a3}{b}
\nccurve[angleA=330,angleB=210,offset=.15]{a3}{a2}
\redmiddlearrow\rput(-.5,-1.2){$1$}
\endpspicture\hspace{.5cm}
\pspicture(-1,-1.5)(1,1)
\pnode(-.8;90){a1}\pnode(-.8;210){a2}\pnode(-.8;330){a3}\pnode(0,0){b}
\ncline{a1}{b}\ncline{a2}{b}\ncline{a3}{b}
\nccurve[angleA=330,angleB=90,offset=-.15]{a3}{a1}
\redmiddlearrow\rput(-.5,-1.2){$1$}
\endpspicture\hspace{.5cm}
\pspicture(-1,-1.5)(1,1)
\pnode(-.8;90){a1}\pnode(-.8;210){a2}\pnode(-.8;330){a3}\pnode(0,0){b}
\ncline{a1}{b}\ncline{a2}{b}\ncline{a3}{b}
\nccurve[angleA=90,angleB=210,offset=-.15]{a1}{a2}
\redmiddlearrow\rput(-.5,-1.2){$1$}
\endpspicture\hspace{.5cm}
\pspicture(-1,-1.5)(1,1)
\pnode(-.8;90){a1}\pnode(-.8;210){a2}\pnode(-.8;330){a3}\pnode(0,0){b}
\ncline{a1}{b}\ncline{a2}{b}\ncline{a3}{b}
\nccurve[angleA=330,angleB=210,offset=.15]{a3}{a2}
\unredmiddlearrow\rput(-.5,-1.2){$v^{-1}$}
\endpspicture\hspace{.5cm}
\pspicture(-1,-1.5)(1,1)
\pnode(-.8;90){a1}\pnode(-.8;210){a2}\pnode(-.8;330){a3}\pnode(0,0){b}
\ncline{a1}{b}\ncline{a2}{b}\ncline{a3}{b}
\nccurve[angleA=330,angleB=90,offset=-.15]{a3}{a1}
\unredmiddlearrow\rput(-.5,-1.2){$v^{-1}$}
\endpspicture\hspace{.5cm}
\pspicture(-1,-1.5)(1,1)
\pnode(-.8;90){a1}\pnode(-.8;210){a2}\pnode(-.8;330){a3}\pnode(0,0){b}
\ncline{a1}{b}\ncline{a2}{b}\ncline{a3}{b}
\nccurve[angleA=90,angleB=210,offset=-.15]{a1}{a2}
\unredmiddlearrow\rput(-.5,-1.2){$v^{-1}$}
\endpspicture
\label{ftweights}
\end{equation}
$$
\pspicture(-1,-1.5)(1,1)
\pnode(-.8;-90){a1}\pnode(-.8;-210){a2}\pnode(-.8;-330){a3}\pnode(0,0){b}
\ncline{a1}{b}\ncline{a2}{b}\ncline{a3}{b}
\nccurve[angleA=-330,angleB=-210,offset=-.15]{a3}{a2}
\redmiddlearrow\rput(-.5,-1.2){$1$}
\endpspicture\hspace{.5cm}
\pspicture(-1,-1.5)(1,1)
\pnode(-.8;-90){a1}\pnode(-.8;-210){a2}\pnode(-.8;-330){a3}\pnode(0,0){b}
\ncline{a1}{b}\ncline{a2}{b}\ncline{a3}{b}
\nccurve[angleA=-330,angleB=-90,offset=.15]{a3}{a1}
\redmiddlearrow\rput(-.5,-1.2){$1$}
\endpspicture\hspace{.5cm}
\pspicture(-1,-1.5)(1,1)
\pnode(-.8;-90){a1}\pnode(-.8;-210){a2}\pnode(-.8;-330){a3}\pnode(0,0){b}
\ncline{a1}{b}\ncline{a2}{b}\ncline{a3}{b}
\nccurve[angleA=-90,angleB=-210,offset=.15]{a1}{a2}
\redmiddlearrow\rput(-.5,-1.2){$1$}
\endpspicture\hspace{.5cm}
\pspicture(-1,-1.5)(1,1)
\pnode(-.8;-90){a1}\pnode(-.8;-210){a2}\pnode(-.8;-330){a3}\pnode(0,0){b}
\ncline{a1}{b}\ncline{a2}{b}\ncline{a3}{b}
\nccurve[angleA=-330,angleB=-210,offset=-.15]{a3}{a2}
\unredmiddlearrow\rput(-.5,-1.2){$v$}
\endpspicture\hspace{.5cm}
\pspicture(-1,-1.5)(1,1)
\pnode(-.8;-90){a1}\pnode(-.8;-210){a2}\pnode(-.8;-330){a3}\pnode(0,0){b}
\ncline{a1}{b}\ncline{a2}{b}\ncline{a3}{b}
\nccurve[angleA=-330,angleB=-90,offset=.15]{a3}{a1}
\unredmiddlearrow\rput(-.5,-1.2){$v$}
\endpspicture\hspace{.5cm}
\pspicture(-1,-1.5)(1,1)
\pnode(-.8;-90){a1}\pnode(-.8;-210){a2}\pnode(-.8;-330){a3}\pnode(0,0){b}
\ncline{a1}{b}\ncline{a2}{b}\ncline{a3}{b}
\nccurve[angleA=-90,angleB=-210,offset=.15]{a1}{a2}
\unredmiddlearrow\rput(-.5,-1.2){$v$}
\endpspicture
$$
For completeness, we give also give a chart of weights
of the $b$'s and $\sigma$'s:
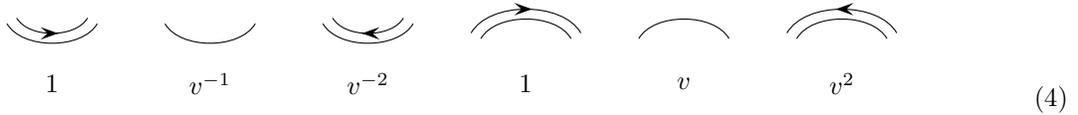
\begin{equation}
\pspicture(-.8,-.5)(.8,1.1)
\pnode(-.6,.6){a1}\pnode(.6,.6){a2}\ncarc[arcangle=-60]{a1}{a2}
\ncarc[arcangle=-60,offset=.15]{a1}{a2}
\redmiddlearrow\rput(0,-.2){$1$}
\endpspicture\hspace{.5cm}
\pspicture(-.8,-.5)(.8,1.1)
\pnode(-.6,.6){a1}\pnode(.6,.6){a2}\ncarc[arcangle=-60]{a1}{a2}
\rput(0,-.2){$v^{-1}$}
\endpspicture\hspace{.5cm}
\pspicture(-.8,-.5)(.8,1.1)
\pnode(-.6,.6){a1}\pnode(.6,.6){a2}\ncarc[arcangle=-60]{a1}{a2}
\ncarc[arcangle=-60,offset=.15]{a1}{a2}
\unredmiddlearrow\rput(0,-.2){$v^{-2}$}
\endpspicture\hspace{.5cm}
\pspicture(-.8,-.5)(.8,1.1)
\pnode(-.6,.4){a1}\pnode(.6,.4){a2}\ncarc[arcangle=60]{a1}{a2}
\ncarc[arcangle=60,offset=.15]{a1}{a2}
\redmiddlearrow\rput(0,-.2){$1$}
\endpspicture\hspace{.5cm}
\pspicture(-.8,-.5)(.8,1.1)\pnode(-.6,.4){a1}\pnode(.6,.4){a2}
\ncarc[arcangle=60]{a1}{a2}
\rput(0,-.2){$v$}
\endpspicture\hspace{.5cm}
\pspicture(-.8,-.5)(.8,1.1)
\pnode(-.6,.4){a1}\pnode(.6,.4){a2}\ncarc[arcangle=60]{a1}{a2}
\ncarc[arcangle=60,offset=.15]{a1}{a2}
\unredmiddlearrow\rput(0,-.2){$v^2$}
\endpspicture
\label{fbweights}
\end{equation}

As an example of computing a state sum using flow lines, the following
are the only two non-zero states with boundary $0,0,0,0,0,1,-1$ in
a certain web $w$:
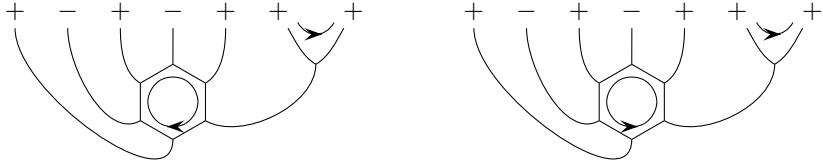
$$
\pspicture(-2.5,-1)(2.6,1.5)
\pnode(.5;270){a1}\pnode(.5;210){a2}\pnode(.5;150){a3}
\pnode(.5; 90){a4}\pnode(.5; 30){a5}\pnode(.5;330){a6}
\pnode(.433;300){d1}\pnode(.433;  0){d2}\pnode(.433; 60){d3}
\pnode(.433;120){d4}\pnode(.433;180){d5}\pnode(.433;240){d6}
\rput(-2.1,1.2){\rnode{b1}{$+$}}\rput(-1.4,1.2){\rnode{b2}{$-$}}
\rput(-.7,1.2){\rnode{b3}{$+$}}\rput(0,1.2){\rnode{b4}{$-$}}
\rput( .7,1.2){\rnode{b5}{$+$}}\pnode(1.9,.5){b6}
\rput(1.4,1.2){\rnode{c1}{$+$}}\rput(2.4,1.2){\rnode{c2}{$+$}}
\ncline{a1}{a2}\ncline{a2}{a3}\ncline{a3}{a4}
\ncline{a4}{a5}\ncline{a5}{a6}\ncline{a6}{a1}
\nccurve[angleA=270,angleB=270,nodesepB=3pt]{a1}{b1}
\nccurve[angleA=210,angleB=270,nodesepB=3pt]{a2}{b2}
\nccurve[angleA=150,angleB=270,nodesepB=3pt]{a3}{b3}
\nccurve[angleA= 90,angleB=270,nodesepB=3pt]{a4}{b4}
\nccurve[angleA= 30,angleB=270,nodesepB=3pt]{a5}{b5}
\nccurve[angleA=330,angleB=270]{a6}{b6}
\nccurve[angleA=150,angleB=300,nodesepB=3pt]{b6}{c1}
\nccurve[angleA= 30,angleB=240,nodesepB=3pt]{b6}{c2}
\nccurve[angleA=300,angleB=240,nodesep=3pt,
    offset=.15,ncurv=1]{c1}{c2}\redmiddlearrow
\nccurve[angleA= 30,angleB=270,offset=.1]{d1}{d2}
\nccurve[angleA= 90,angleB=330,offset=.1]{d2}{d3}
\nccurve[angleA=150,angleB= 30,offset=.1]{d3}{d4}
\nccurve[angleA=210,angleB= 90,offset=.1]{d4}{d5}
\nccurve[angleA=270,angleB=150,offset=.1]{d5}{d6}
\nccurve[angleA=330,angleB=210,offset=.1]{d6}{d1}
\unredmiddlearrow
\endpspicture
\hspace{1cm}
\pspicture(-2.5,-1)(2.6,1.5)
\pnode(.5;270){a1}\pnode(.5;210){a2}\pnode(.5;150){a3}
\pnode(.5; 90){a4}\pnode(.5; 30){a5}\pnode(.5;330){a6}
\pnode(.433;300){d1}\pnode(.433;  0){d2}\pnode(.433; 60){d3}
\pnode(.433;120){d4}\pnode(.433;180){d5}\pnode(.433;240){d6}
\rput(-2.1,1.2){\rnode{b1}{$+$}}\rput(-1.4,1.2){\rnode{b2}{$-$}}
\rput(-.7,1.2){\rnode{b3}{$+$}}\rput(0,1.2){\rnode{b4}{$-$}}
\rput( .7,1.2){\rnode{b5}{$+$}}\pnode(1.9,.5){b6}
\rput(1.4,1.2){\rnode{c1}{$+$}}\rput(2.4,1.2){\rnode{c2}{$+$}}
\ncline{a1}{a2}\ncline{a2}{a3}\ncline{a3}{a4}
\ncline{a4}{a5}\ncline{a5}{a6}\ncline{a6}{a1}
\nccurve[angleA=270,angleB=270,nodesepB=3pt]{a1}{b1}
\nccurve[angleA=210,angleB=270,nodesepB=3pt]{a2}{b2}
\nccurve[angleA=150,angleB=270,nodesepB=3pt]{a3}{b3}
\nccurve[angleA= 90,angleB=270,nodesepB=3pt]{a4}{b4}
\nccurve[angleA= 30,angleB=270,nodesepB=3pt]{a5}{b5}
\nccurve[angleA=330,angleB=270]{a6}{b6}
\nccurve[angleA=150,angleB=300,nodesepB=3pt]{b6}{c1}
\nccurve[angleA= 30,angleB=240,nodesepB=3pt]{b6}{c2}
\nccurve[angleA=300,angleB=240,nodesep=3pt,
    offset=.15,ncurv=1]{c1}{c2}\redmiddlearrow
\nccurve[angleA= 30,angleB=270,offset=.1]{d1}{d2}
\nccurve[angleA= 90,angleB=330,offset=.1]{d2}{d3}
\nccurve[angleA=150,angleB= 30,offset=.1]{d3}{d4}
\nccurve[angleA=210,angleB= 90,offset=.1]{d4}{d5}
\nccurve[angleA=270,angleB=150,offset=.1]{d5}{d6}
\nccurve[angleA=330,angleB=210,offset=.1]{d6}{d1}
\redmiddlearrow
\endpspicture
$$
Since the weights of these states are $v^{-1}$ and $v^{-3}$, we therefore
conclude that the coefficient of
$$e^+_0 \tensor e^-_0 \tensor e^+_0 \tensor e^-_0
\tensor e^+_0 \tensor e^+_1 \tensor e^+_{-1}$$
in $w$ is $v^{-1} + v^{-3}$ in the tensor product basis.

Finally, note that the weight of any state of any monomial web
is either a power of $v$ or zero.  Thus, weights cannot cancel
in state sums, and any state sum takes values in $\N[v,v^{-1}]$.

\section{From paths and strings to non-elliptic webs}

In order to compare the web and dual canonical bases, we must
enumerate non-elliptic webs by certain strings of elements of
$\{-1,0,1\}$, namely those that correspond to weight lattice paths
confined to a Weyl chamber of $\sl(3)$.

More precisely, let $S = s_1,\ldots,s_n$ be a string of signs, and let
$J = j_1,\ldots,j_n$ be a string of states chosen from $\{-1,0,1\}$.  Each
vector $e^{s_k}_{j_k}$ has a weight $\mu_k$, and we may define a path $0 =
\pi_0,\pi_1,\pi_2,\ldots,\pi_n$ in the weight lattice of $\sl(3)$ such that
$\pi_k = \mu_k + \pi_{k-1}$. The dominant Weyl chamber is defined as the subset
of the weight lattice consisting of positive integral linear combinations of
the weights $\mu^+$ and $\mu^-$ of $e^+_1$ and $e^-_1$. It is well-known
that, for fixed $s_1,\ldots,s_n$, the number of strings $j_1,\ldots,j_n$ that
produce a lattice path in the dominant Weyl chamber from the origin to itself
equals
$$\dim \Inv(V^{s_1} \tensor V^{s_2} \tensor \ldots \tensor V^{s_n}).$$
Call such lattice paths dominant.

Given a sign string and a string of states, we define a web by inductive
rules called the growth algorithm:
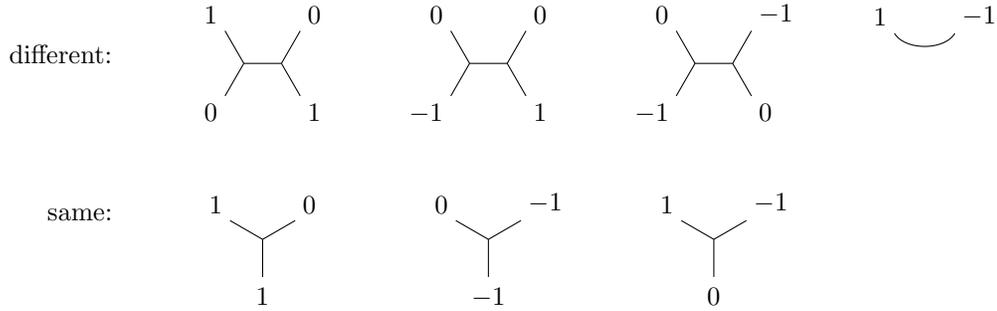
\begin{align*}
\mbox{different:} & \hspace{1cm}
\pspicture[.5](-1,-.7)(1,.7)
\qline( .25,0)( .5,.433)\qline( .25,0)( .5,-.433)
\qline(-.25,0)(-.5,.433)\qline(-.25,0)(-.5,-.433)
\qline(-.25,0)(.25,0)
\rput[br](-.6, .533){$1$}\rput[bl](.6, .533){$0$}
\rput[tr](-.6,-.533){$0$}\rput[tl](.6,-.533){$1$}
\endpspicture\hspace{1cm}
\pspicture[.5](-1,-.7)(1,.7)
\qline( .25,0)( .5,.433)\qline( .25,0)( .5,-.433)
\qline(-.25,0)(-.5,.433)\qline(-.25,0)(-.5,-.433)
\qline(-.25,0)(.25,0)
\rput[br](-.6, .533){$0$}\rput[bl](.6, .533){$0$}
\rput[tr](-.6,-.533){$-1$}\rput[tl](.6,-.533){$1$}
\endpspicture\hspace{1cm}
\pspicture[.5](-1,-.7)(1,.7)
\qline( .25,0)( .5,.433)\qline( .25,0)( .5,-.433)
\qline(-.25,0)(-.5,.433)\qline(-.25,0)(-.5,-.433)
\qline(-.25,0)(.25,0)
\rput[br](-.6, .533){$0$}\rput[bl](.6, .533){$-1$}
\rput[tr](-.6,-.533){$-1$}\rput[tl](.6,-.533){$0$}
\endpspicture\hspace{1cm}
\pspicture[.5](-.8,-.7)(.8,.7)
\pccurve[angleA=-60,angleB=240](-.4,.4)(.4,.4)
\rput[br](-.5,.5){$1$}\rput[bl](.5,.5){$-1$}
\endpspicture
\\
\mbox{same:} & \hspace{1cm}
\pspicture[.5](-1,-1)(1,1.5)
\pnode(-.5;90){a1}\pnode(-.5;210){a2}\pnode(-.5;330){a3}\pnode(0,0){b}
\ncline{a1}{b}\ncline{a2}{b}\ncline{a3}{b}
\rput[br](-.533,.35){$1$}\rput[bl](.533,.35){$0$}
\rput[t](0,-.65){$1$}
\endpspicture\hspace{1cm}
\pspicture[.5](-1,-1)(1,1.5)
\pnode(-.5;90){a1}\pnode(-.5;210){a2}\pnode(-.5;330){a3}\pnode(0,0){b}
\ncline{a1}{b}\ncline{a2}{b}\ncline{a3}{b}
\rput[br](-.533,.35){$0$}\rput[bl](.533,.35){$-1$}
\rput[t](0,-.65){$-1$}
\endpspicture\hspace{1cm}
\pspicture[.5](-1,-1)(1,1.5)
\pnode(-.5;90){a1}\pnode(-.5;210){a2}\pnode(-.5;330){a3}\pnode(0,0){b}
\ncline{a1}{b}\ncline{a2}{b}\ncline{a3}{b}
\rput[br](-.533,.35){$1$}\rput[bl](.533,.35){$-1$}
\rput[t](0,-.65){$0$}
\endpspicture
\end{align*}
Initially, the web consists of parallel strands whose orientations are
given by the sign string. The rules indicate that if the state string
admits a substring of the type appearing at the top in one of the cases
(taken from the top row if the two signs are different and from the
bottom row if they are the same), we should concatenate the
corresponding web and replace the substrings with what appears at the
bottom of the web.  If none of the patterns at the top appear anywhere,
the growth algorithm stops.  For example, the growth algorithm converts
the sign string $+-+-+++$ and the state string $1,1,0,0,-1,0,-1$ to the
web:
$$
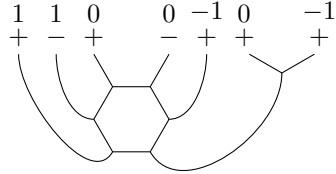
\pspicture(-1.9,-1)(2.9,1.7)
\psline(.5;0)(.5;60)(.5;120)(.5;180)(.5;240)(.5;300)(.5;0)
\psline(.5;60)(1;60)
\psline(.5;120)(1;120)
\pccurve[angleA=180,angleB=270](.5;180)(-1,.866)
\pccurve[angleA=0,angleB=270](.5;0)(1,.866)
\pccurve[angleA=240,angleB=270](.5;240)(-1.5,.866)
\pccurve[angleA=300,angleB=270,ncurv=.8](.5;300)(2,.616)
\psline(1.567,.866)(2,.616)(2.433,.866)
\rput[b](-1.5,.966){$+$}\rput[b](-1.5,1.3){$1$}
\rput[b](-1,.966){$-$}\rput[b](-1,1.3){$1$}
\rput[b](-.5,.966){$+$}\rput[b](-.5,1.3){$0$}
\rput[b](.5,.966){$-$}\rput[b](.5,1.3){$0$}
\rput[b](1,.966){$+$}\rput[b](1,1.3){$-1$}
\rput[b](1.5,.966){$+$}\rput[b](1.5,1.3){$0$}
\rput[b](2.5,.966){$+$}\rput[b](2.5,1.3){$-1$}
\endpspicture
$$
In this case the growth algorithm continues until the sign and state
strings have length 0.

The validity of the growth algorithm rests on the following lemmas.

\begin{lemma} Given any sign and state string, the web produced
by the growth algorithm does not depend on the order in which the
substrings are replaced.
\end{lemma}
\begin{proof} (Sketch) The proof is by induction.  Order state strings
by their length; if two state strings have the same length, order them
lexicographically.  A minimal counterexample consisting of a sign
string $S$ and a state string $J$ must have two different replaceable
substrings that ultimately result in two different webs; given two such
replacements $r_1$ and $r_2$, let $w_1$ and $w_2$ be the two webs that
result, and let $S_1$ and $S_2$ and $J_1$ and $J_2$ be the sign and
state strings that result.  The pairs $(S_1,J_1)$ and $(S_2,J_2)$ are
not counterexamples because $J_1$ and $J_2$ come before $J$; therefore
the growth algorithm is order-independent for both of these strings. To
obtain a contradiction, it suffices to complete a diamond by finding a
pair $(S_3,J_3)$ which can be obtained from either $(S_1,J_1)$ or
$(S_2,J_2)$ by the growth algorithm, for example:
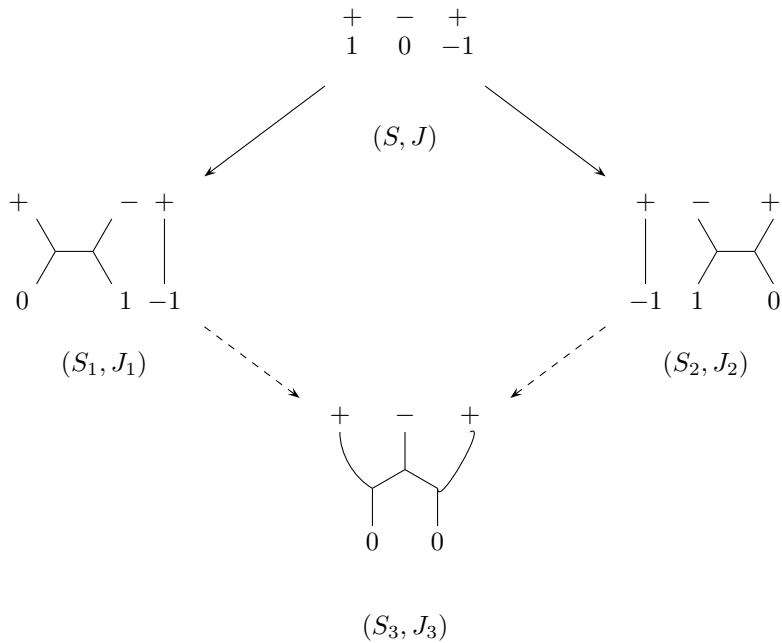
\begin{equation}
\pspicture(-5.5,-6)(5.5,4.5)
\rput(0,3){\rnode{a1}{\pspicture(-1,-.5)(1,.5)
\rput[b](-.7,0){$+$}\rput[b](0,0){$-$}\rput[b](.7,0){$+$}
\rput[t](-.7,-.15){$1$}\rput[t](0,-.15){$0$}\rput[t](.7,-.15){$-1$}
\endpspicture}}\rput(0,1.5){$(S,J)$}
\rput(-4,0){\rnode{a2}{\pspicture(-.7,-.7)(1.5,.7)
\qline( .25,0)( .5,.433)\qline( .25,0)( .5,-.433)
\qline(-.25,0)(-.5,.433)\qline(-.25,0)(-.5,-.433)\qline(-.25,0)(.25,0)
\qline(1.2,.433)(1.2,-.433)
\rput[br](-.6, .533){$+$}\rput[bl](.6, .533){$-$}\rput[b](1.2,.533){$+$}
\rput[tr](-.6,-.533){$0$}\rput[tl](.6,-.533){$1$}\rput[t](1.2,-.533){$-1$}
\endpspicture}}\rput(-4,-1.5){$(S_1,J_1)$}
\rput(4,0){\rnode{a3}{\pspicture(-1.5,-.7)(.7,.7)
\qline(-.25,0)(-.5,.433)\qline(-.25,0)(-.5,-.433)
\qline(.25,0)(.5,.433)\qline(.25,0)(.5,-.433)\qline(.25,0)(-.25,0)
\qline(-1.2,.433)(-1.2,-.433)
\rput[br](.6, .533){$+$}\rput[bl](-.6, .533){$-$}\rput[b](-1.2,.533){$+$}
\rput[tr](.6,-.533){$0$}\rput[tl](-.6,-.533){$1$}\rput[t](-1.2,-.533){$-1$}
\endpspicture}}\rput(4,-1.5){$(S_2,J_2)$}
\rput(0,-3){\rnode{a4}{\pspicture(-1.1,-1.1)(1.1,.9)
\pccurve[angleA=270,angleB=150](-.866,.5)(-.433,-.25)
\pccurve[angleA=270,angleB=30](.433,-.25)(.866,.5)
\psline(-.433,-.25)(0,0)(.433,-.25)
\psline(0,0)(0,.5)\psline(-.433,-.25)(-.433,-.75)
\psline(.433,-.25)(.433,-.75)
\rput[b](-.866,.6){$+$}\rput[b](0,.6){$-$}\rput[b](.866,.6){$+$}
\rput[t](-.433,-.85){$0$}\rput[t](.433,-.85){$0$}
\endpspicture}}\rput(0,-5){$(S_3,J_3)$}
\ncline[arrows=->,nodesep=.3]{a1}{a2}
\ncline[arrows=->,nodesep=.3]{a1}{a3}
\ncline[arrows=->,nodesep=.3,linestyle=dashed]{a2}{a4}
\ncline[arrows=->,nodesep=.3,linestyle=dashed]{a3}{a4}
\endpspicture\label{fdiamond}
\end{equation}
If the replacements $r_1$ and $r_2$ have disjoint substrings, then we can
trivially complete the diamond by appling $r_2$ after $r_1$ and vice-versa.
There is a short list of cases in which they are not disjoint, and we can
complete the diamond on a case-by-case basis.  Figure~(\ref{fdiamond}) gives
one of the cases.
\end{proof}

\begin{lemma} Any web produced by the growth algorithm is non-elliptic.
\end{lemma}

\begin{proof}(Sketch) In the growth rules, an internal face can only be
``born'' with a rule that produces an H.  The only way to obtain a
2-sided face would be to close off the face immediately with a U.
However, the indices that result from attaching an H rule out this
possibility.

A square face can be ruled out by a more complicated version of the same
reasoning.  A square might hypothetically have one of four possible
histories:
$$
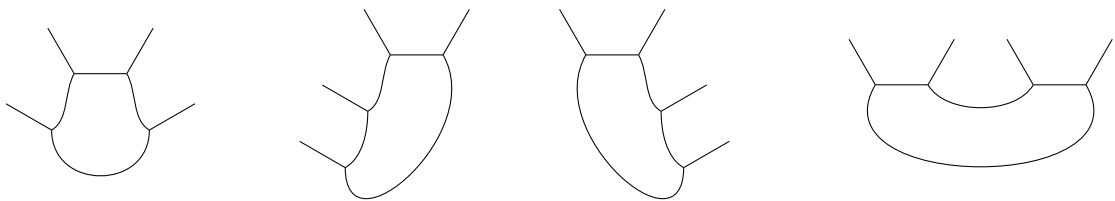
\pspicture[.5](-1.3,-1.6)(1.3,.8)
\pnode(-.35,0){a1}\pnode(-.65,-.75){a2}
\pnode(.65,-.75){a3}\pnode(.35,0){a4}
\nccurve[angleA=240,angleB=30]{a1}{a2}
\nccurve[angleA=270,angleB=270,ncurv=1.25]{a2}{a3}
\nccurve[angleA=150,angleB=300]{a3}{a4}\ncline{a4}{a1}
\psline(a1)([nodesep=.7,angle=120]a1)\psline(a2)([nodesep=.7,angle=150]a2)
\psline(a3)([nodesep=.7,angle=30]a3)\psline(a4)([nodesep=.7,angle=60]a4)
\endpspicture
\hspace{1cm}
\pspicture[.5](-1.9,-2.1)(.8,.8)
\pnode(-.35,0){a1}\pnode(-.65,-.75){a2}
\pnode(-.95,-1.5){a3}\pnode(.35,0){a4}
\nccurve[angleA=240,angleB=30]{a1}{a2}
\nccurve[angleA=270,angleB=30]{a2}{a3}
\nccurve[angleA=270,angleB=300,ncurv=1.2]{a3}{a4}\ncline{a4}{a1}
\psline(a1)([nodesep=.7,angle=120]a1)\psline(a2)([nodesep=.7,angle=150]a2)
\psline(a3)([nodesep=.7,angle=150]a3)\psline(a4)([nodesep=.7,angle=60]a4)
\endpspicture\hspace{1cm}
\pspicture[.5](-.8,-2.1)(1.9,.8)
\pnode(.35,0){a1}\pnode(.65,-.75){a2}
\pnode(.95,-1.5){a3}\pnode(-.35,0){a4}
\nccurve[angleA=300,angleB=150]{a1}{a2}
\nccurve[angleA=270,angleB=150]{a2}{a3}
\nccurve[angleA=270,angleB=240,ncurv=1.2]{a3}{a4}\ncline{a4}{a1}
\psline(a1)([nodesep=.7,angle=60]a1)\psline(a2)([nodesep=.7,angle=30]a2)
\psline(a3)([nodesep=.7,angle=30]a3)\psline(a4)([nodesep=.7,angle=120]a4)
\endpspicture\hspace{1cm}
\pspicture[.5](-2,-1.3)(2,.8)
\pnode(-1.4,0){a1}\pnode(-.7,0){a2}
\pnode(.7,0){a3}\pnode(1.4,0){a4}
\ncline{a1}{a2}\ncline{a3}{a4}
\nccurve[angleA=300,angleB=240]{a2}{a3}
\nccurve[angleA=240,angleB=300,ncurv=1.2]{a1}{a4}
\psline(a1)([nodesep=.7,angle=120]a1)\psline(a2)([nodesep=.7,angle=60]a2)
\psline(a3)([nodesep=.7,angle=120]a3)\psline(a4)([nodesep=.7,angle=60]a4)
\endpspicture
$$
In the first three cases, the extra vertices may belong either to H's or
Y's produced by the growth algorithm at adjacent locations. Working
backwards from the final U, one quickly concludes that none of the
histories are possible.
\end{proof}

\begin{lemma} If a sign and state string correspond to a
dominant lattice path, then the growth algorithm does not terminate
until the strings have length 0. \label{lterm}
\end{lemma}
\begin{proof} (Sketch) The proof is again by induction on length and
lexicographic order.  Observe that the growth algorithm only terminates
at a non-decreasing state string (and only then when all positions $k$
such that $j_k = 0$ have the same sign $s_k$).  On the other hand, the
state string of a dominant path of length greater than 0 must begin
with 1 and end with -1.  Thus, it suffices to show that a growth rule
applied to a dominant path produces another dominant path.  None of the
14 growth rules change the endpoints of the corresponding lattice path,
and each of them either reduces the set of vertices it visits or
modifies it in a way that cannot lead to an excursion outside of the
dominant Weyl chamber.  For example, a growth rule that produces a Y
replaces two consecutive steps of the path by one step:
$$
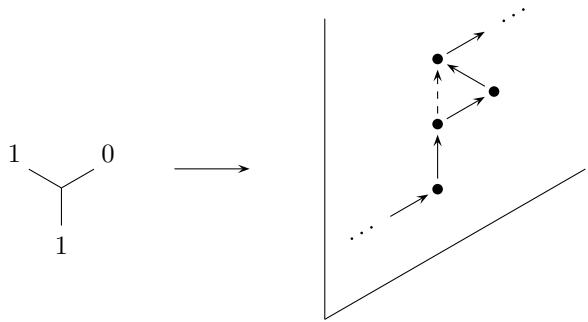
\pspicture[.5](-1,-1)(1,1.5)
\pnode(-.5;90){a1}\pnode(-.5;210){a2}\pnode(-.5;330){a3}\pnode(0,0){b}
\ncline{a1}{b}\ncline{a2}{b}\ncline{a3}{b}
\rput[br](-.533,.35){$1$}\rput[bl](.533,.35){$0$}
\rput[t](0,-.65){$1$}
\endpspicture
\psgoesto
\pspicture[.5](-.5,0)(4,4)
\psline(0,0)(0,4)
\psline(0,0)(4;30)
\pnode(.75,1.299){a1}
\pnode(1.5,1.732){a2}
\pnode(1.5,2.598){a3}
\pnode(2.25,3.031){a4}
\pnode(1.5,3.464){a5}
\pnode(2.25,3.897){a6}
\ncline[arrows=->,nodesep=4pt]{a1}{a2}
\ncline[arrows=->,nodesep=4pt]{a2}{a3}
\ncline[arrows=->,nodesep=4pt]{a3}{a4}
\ncline[arrows=->,nodesep=4pt]{a4}{a5}
\ncline[arrows=->,nodesep=4pt]{a5}{a6}
\ncline[arrows=->,linestyle=dashed,nodesep=4pt]{a3}{a5}
\pscircle*(a2){2pt}
\pscircle*(a3){2pt}
\pscircle*(a4){2pt}
\pscircle*(a5){2pt}
\pscircle*([nodesep=.15,angle=210]a1){.5pt}
\pscircle*([nodesep=.3,angle=210]a1){.5pt}
\pscircle*([nodesep=.45,angle=210]a1){.5pt}
\pscircle*([nodesep=.15,angle=30]a6){.5pt}
\pscircle*([nodesep=.3,angle=30]a6){.5pt}
\pscircle*([nodesep=.45,angle=30]a6){.5pt}
\endpspicture
$$
\end{proof}

The converse of Lemma~\ref{lterm} also holds \cite{Khovanov}.

The growth algorithm has a notable inverse for dominant paths. Let $w$ be a
non-elliptic web.  Given points $P$ and $Q$ on the boundary of $w$ lying
between endpoints, a minimal cut path is a transversely oriented arc from $P$
to $Q$ which is transverse to $w$ and which crosses as few strands as
possible:
$$
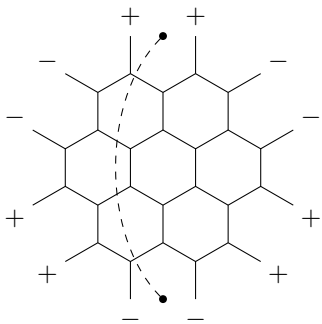
\pspicture[.47](-2.2,-2.2)(2.2,2.2)
\psline(-1.299,-1.25)(-.866,-1)(-.433,-1.25)
    (0,-1)(.433,-1.25)(.866,-1)(1.299,-1.25)
\psline(-1.732,-.5)(-1.299,-.25)(-.866,-.5)(-.433,-.25)(0,-.5)(.433,-.25)
    (.866,-.5)(1.299,-.25)(1.732,-.5)
\psline(-1.732,.5)(-1.299,.25)(-.866,.5)(-.433,.25)
    (0,.5)(.433,.25)(.866,.5)(1.299,.25)(1.732,.5)
\psline(-1.299,1.25)(-.866,1)(-.433,1.25)(0,1)(.433,1.25)(.866,1)(1.299,1.25)
\psline(-.433,1.75)(-.433,1.25)\psline(.433,1.75)(.433,1.25)
\psline(-.866,1)(-.866,.5)\psline(0,1)(0,.5)\psline(.866,1)(.866,.5)
\psline(-1.299,.25)(-1.299,-.25)\psline(-.433,.25)(-.433,-.25)
\psline(1.299,.25)(1.299,-.25)\psline(.433,.25)(.433,-.25)
\psline(-.866,-1)(-.866,-.5)\psline(0,-1)(0,-.5)\psline(.866,-1)(.866,-.5)
\psline(-.433,-1.75)(-.433,-1.25)\psline(.433,-1.75)(.433,-1.25)
\rput[b](-.433,1.9){$+$}\rput[b](.433,1.9){$+$}
\rput[br](-1.399,1.3){$-$}\rput[bl](1.399,1.3){$-$}
\rput[br](-1.832,.55){$-$}\rput[bl](1.832,.55){$-$}
\rput[tr](-1.832,-.55){$+$}\rput[tl](1.832,-.55){$+$}
\rput[tr](-1.399,-1.3){$+$}\rput[tl](1.399,-1.3){$+$}
\rput[t](-.433,-1.9){$-$}\rput[t](.433,-1.9){$-$}
\pcarc[arrows=*-*,linestyle=dashed,arcangle=45](0,-1.75)(0,1.75)
\endpspicture
$$
The weight of a minimal cut path is $a\mu^+ + b\mu^-$ if the cut path
crosses $a$ strands whose orientations agree with that of the arc and
$b$ strands whose orientations disagree. Although minimal cut paths are
not necessarily unique, their weights are \cite{Kuperberg:spiders}.
Moreover, any two minimal cut paths are connected by a sequence of
H-moves:
$$
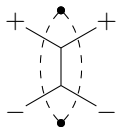
\pspicture(-1,-.9)(1,.9)
\pnode(0,.25){a1}\pnode(0,-.25){a2}
\rput([nodesep=.7,angle= 30]a1){\rnode{b1}{$+$}}
\rput([nodesep=.7,angle=150]a1){\rnode{b2}{$+$}}
\rput([nodesep=.7,angle=330]a2){\rnode{b3}{$-$}}
\rput([nodesep=.7,angle=210]a2){\rnode{b4}{$-$}}
\ncline{a1}{b1}\ncline{a1}{b2}\ncline{a1}{a2}
\ncline{a2}{b3}\ncline{a2}{b4}
\pcarc[arrows=*-*,linestyle=dashed,arcangle=45]
    (0,.75)(0,-.75)
\pcarc[arrows=*-*,linestyle=dashed,arcangle=-45]
    (0,.75)(0,-.75)
\endpspicture
$$
Now let $w$ be a non-elliptic web with $n$ endpoints that are linearly
ordered (rather than cyclically ordered) and lie above $w$, as might be
produced by the growth algorithm. Let $P$ be a point below $w$, and let
$Q_0,Q_1,\ldots,Q_n$ be points that alternate with the endpoints of
$w$.  Let $\pi_k$ be the weight of a minimal cut path from $P$ to
$Q_k$:
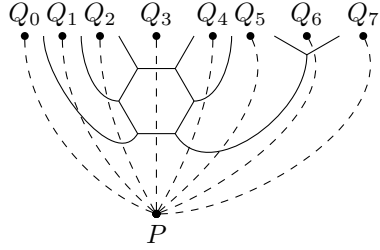
$$
\pspicture(-1.9,-2)(2.9,1.7)
\psline(.5;0)(.5;60)(.5;120)(.5;180)(.5;240)(.5;300)(.5;0)
\psline(.5;60)(1;60)
\psline(.5;120)(1;120)
\pccurve[angleA=180,angleB=270](.5;180)(-1,.866)
\pccurve[angleA=0,angleB=270](.5;0)(1,.866)
\pccurve[angleA=240,angleB=270](.5;240)(-1.5,.866)
\pccurve[angleA=300,angleB=270,ncurv=.8](.5;300)(2,.616)
\psline(1.567,.866)(2,.616)(2.433,.866)
\pnode(0,-1.5){p}\rput[t](0,-1.65){$P$}
\pccurve[arrows=*-*,linestyle=dashed,%
    angleA=160,angleB=270](p)(-1.75,.866)\rput[b](-1.75,1){$Q_0$}
\pccurve[arrows=*-*,linestyle=dashed,%
    angleA=140,angleB=270](p)(-1.25,.866)\rput[b](-1.25,1){$Q_1$}
\pccurve[arrows=*-*,linestyle=dashed,%
    angleA=120,angleB=270](p)(-.75,.866)\rput[b](-.75,1){$Q_2$}
\pccurve[arrows=*-*,linestyle=dashed,%
    angleA=90,angleB=270](p)(0,.866)\rput[b](0,1){$Q_3$}
\pccurve[arrows=*-*,linestyle=dashed,%
    angleA=60,angleB=270](p)(.75,.866)\rput[b](.75,1){$Q_4$}
\pccurve[arrows=*-*,linestyle=dashed,%
    angleA=40,angleB=295](p)(1.25,.866)\rput[b](1.25,1){$Q_5$}
\pccurve[arrows=*-*,linestyle=dashed,%
    angleA=20,angleB=300](p)(2,.866)\rput[b](2,1){$Q_6$}
\pccurve[arrows=*-*,linestyle=dashed,%
    angleA=0,angleB=295](p)(2.75,.866)\rput[b](2.75,1){$Q_7$}
\endpspicture
$$
The sequence $0 = \pi_0,\pi_1,\ldots,\pi_n = 0$ is a dominant
lattice path \cite{Kuperberg:spiders}.  Each successive difference
$\pi_k - \pi_{k-1}$ is the weight of some $e^{s_k}_{j_k}$, so we can
reconstruct both a sign string (which is given directly by the boundary
of $w$) and a state string from the web $w$. Call this procedure the
minimal cut path algorithm.

\begin{proposition} The minimal cut path and growth algorithms are inverses.
\label{pinverse}
\end{proposition}
\begin{proof}
For each fixed sign string $S$, let $w^S_J$ be the web produced by the
growth algorithm from the state string $J = j_1,\ldots,j_n$, and let
$m(w)$ be the state string produced by the minimal cut path algorithm
from the basis web $w$.  We will show that $m \circ g = I$.
Since $m$ is a bijection \cite{Kuperberg:spiders}, it follows
that $g \circ m = I$ also.

We extend the growth rules to create a system
of flow lines and minimal cut paths along with the basis web.
The extension is given by the following diagrams:
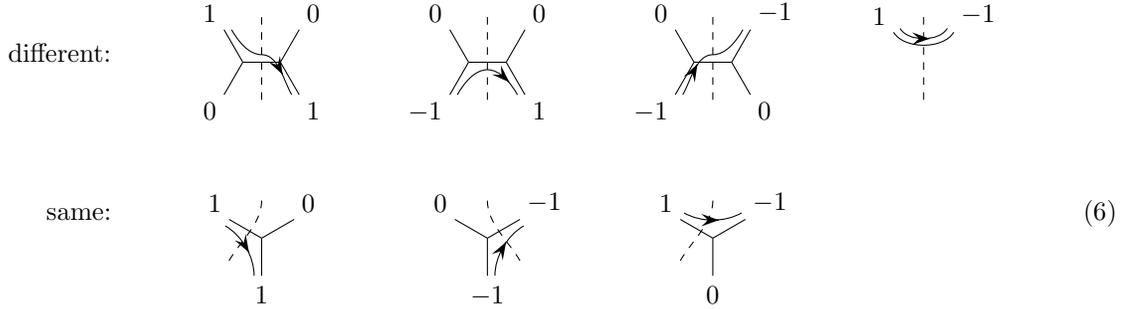
\begin{align}
\mbox{different:} & \hspace{1cm}
\pspicture[.5](-1,-.7)(1,.7)
\qline( .25,0)( .5,.433)\qline( .25,0)( .5,-.433)
\qline(-.25,0)(-.5,.433)\qline(-.25,0)(-.5,-.433)\qline(-.25,0)(.25,0)
\rput[br](-.6, .533){$1$}\rput[bl](.6, .533){$0$}
\rput[tr](-.6,-.533){$0$}\rput[tl](.6,-.533){$1$}
\pccurve[angleA=300,angleB=180](-.4,.433)(0,.1)
\pccurve[angleA=0,angleB=120](0,.1)(.4,-.433)
\redmiddlearrow
\psline[linestyle=dashed](0,.6)(0,-.5)
\endpspicture\hspace{1cm}
\pspicture[.5](-1,-.7)(1,.7)
\qline( .25,0)( .5,.433)\qline( .25,0)( .5,-.433)
\qline(-.25,0)(-.5,.433)\qline(-.25,0)(-.5,-.433)\qline(-.25,0)(.25,0)
\rput[br](-.6, .533){$0$}\rput[bl](.6, .533){$0$}
\rput[tr](-.6,-.533){$-1$}\rput[tl](.6,-.533){$1$}
\pccurve[angleA=60,angleB=180](-.4,-.433)(0,-.1)
\pccurve[angleA=0,angleB=120](0,-.1)(.4,-.433)
\redmiddlearrow
\psline[linestyle=dashed](0,.6)(0,-.5)
\endpspicture\hspace{1cm}
\pspicture[.5](-1,-.7)(1,.7)
\qline( .25,0)( .5,.433)\qline( .25,0)( .5,-.433)
\qline(-.25,0)(-.5,.433)\qline(-.25,0)(-.5,-.433)\qline(-.25,0)(.25,0)
\rput[br](-.6, .533){$0$}\rput[bl](.6, .533){$-1$}
\rput[tr](-.6,-.533){$-1$}\rput[tl](.6,-.533){$0$}
\pccurve[angleA=240,angleB=0](.4,.433)(0,.1)
\pccurve[angleA=180,angleB=60](0,.1)(-.4,-.433)
\unredmiddlearrow
\psline[linestyle=dashed](0,.6)(0,-.5)
\endpspicture\hspace{1cm}
\pspicture[.5](-.8,-.7)(.8,.7)
\pnode(-.4,.4){a1}\pnode(.4,.4){a2}
\nccurve[angleA=-60,angleB=240]{a1}{a2}
\rput[br](-.5,.5){$1$}\rput[bl](.5,.5){$-1$}
\ncarc[arcangle=-60,offset=.1]{a1}{a2}
\redmiddlearrow
\psline[linestyle=dashed](0,.6)(0,-.5)
\endpspicture \nonumber \\
\mbox{same:} & \hspace{1cm}
\pspicture[.5](-1,-1)(1,1.5)
\pnode(-.5;90){a1}\pnode(-.5;210){a2}\pnode(-.5;330){a3}\pnode(0,0){b}
\ncline{a1}{b}\ncline{a2}{b}\ncline{a3}{b}
\rput[br](-.533,.35){$1$}\rput[bl](.533,.35){$0$}
\rput[t](0,-.65){$1$}
\nccurve[angleA=330,angleB=90,offset=-.1]{a3}{a1}
\redmiddlearrow
\pccurve[angleA=270,angleB=60,linestyle=dashed]%
    (0,.5)(-.433,-.3)
\endpspicture\hspace{1cm}
\pspicture[.5](-1,-1)(1,1.5)
\pnode(-.5;90){a1}\pnode(-.5;210){a2}\pnode(-.5;330){a3}\pnode(0,0){b}
\ncline{a1}{b}\ncline{a2}{b}\ncline{a3}{b}
\rput[br](-.533,.35){$0$}\rput[bl](.533,.35){$-1$}
\rput[t](0,-.65){$-1$}
\nccurve[angleA=90,angleB=210,offset=-.1]{a1}{a2}
\redmiddlearrow
\pccurve[angleA=270,angleB=120,linestyle=dashed]%
    (0,.5)(.433,-.3)
\endpspicture\hspace{1cm}
\pspicture[.5](-1,-1)(1,1.5)
\pnode(-.5;90){a1}\pnode(-.5;210){a2}\pnode(-.5;330){a3}\pnode(0,0){b}
\ncline{a1}{b}\ncline{a2}{b}\ncline{a3}{b}
\rput[br](-.533,.35){$1$}\rput[bl](.533,.35){$-1$}
\rput[t](0,-.65){$0$}
\nccurve[angleA=330,angleB=210,offset=.1]{a3}{a2}
\redmiddlearrow
\pccurve[angleA=270,angleB=60,linestyle=dashed]%
    (0,.5)(-.433,-.3)
\endpspicture \label{fextended}
\end{align}
In each case, cut paths may merge, so that at any step after the first one, a
single cut path represented in the diagram may be replaced by many parallel
cut paths. In the final result, the cut paths are all minimal by the
principle that $\mbox{MIN CUT} \ge \mbox{MAX FLOW}$. \Ie, if the sum of
the indices to the left of a given cut path $\gamma$ is $m$, then by
conservation of flow, $\gamma$ must cross at least $m$ flow lines.  By
construction, any strand that $\gamma$ cross has flow on it and all flow
across $\gamma$ is to the right; therefore $\gamma$ must cross exactly $m$
flow lines.

Since the chosen cut paths are minimal, it is routine to show by induction on
the number of steps in the growth algorithm that the minimal cut path
algorithm reconstructs the original state string.
\end{proof}

\begin{proof}[Remark] A well-known theorem in the theory of linear
programming is often summarized by the maxim $\mbox{MIN CUT} =
\mbox{MAX FLOW}$. However, if we applied this theorem directly, we
would have to allow the possibility of fractional flow lines.  Thus, we
have shown in our case that the maximal flow in the linear sense can be
achieved combinatorially without using flow lines with forks.
\end{proof}

Reference~\citen{Kuperberg:spiders} demonstrated that the set of non-elliptic
webs is a $\C(q)$-basis of $\Inv(V^{s_1} \tensor \ldots \tensor V^{s_n})$.
Proposition~\ref{pinverse} gives a way to index this basis by dominant paths. 
Here we can obtain a stronger result. If $S = s_1,\ldots,s_n$ is a sign
string and $J = j_1,\ldots,j_n$ is a state string, define $e^S_J$ by:
$$e^S_J = e^{s_1}_{j_1} \tensor e^{s_2}_{j_2} \tensor \ldots
\tensor e^{s_n}_{j_n}$$

\begin{theorem} The tensor $w^S_J$ expands as
$$w^S_J = e^S_J + \sum_{J' < J} c(S,J,J') e^S_{J'}$$
for some coefficients $c(S,J,J') \in \N[v,v^{-1}]$,
where the state strings $J$ and $J'$ are ordered lexicographically.
\label{thleading}
\end{theorem}

It follows that the non-elliptic webs are a $\Z[v,v^{-1}]$-basis of the
invariant spaces in which they live.

\begin{proof}
The result follows from the existence of the minimal cut paths. Recall
that the coefficients $c(S,J,J')$ are state sums. For each $k$, the cut
path $\gamma$ from $P$ to $Q_k$ cuts $m$ strands, where $m = \sum_{\ell
= 1}^k j_\ell$.  But if $J' > J$, we can choose the first $k$ such that
$j_k \ne j'_k$; in this case, $j'_k > j_k$.  The cut path $\gamma$ must
cut at least $\sum_{\ell = 1}^k j'_\ell > m$ flow lines in any state
contributing to $c(S,J,J')$. But this is impossible, since $\gamma$
only cuts $m$ strands.

The case $J' = J$ is more delicate.  We claim that only one non-zero
state contributes to $c(S,J,J)$, namely the state $x$ constructed in
Proposition~\ref{pinverse}.  Since $\mbox{MIN CUT} = \mbox{MAX FLOW}$
by this proposition, every edge which intersects any minimal cut path
must carry flow in every non-zero state.  The claim follows if we can
show that minimal cut paths meet every edge that carries flow in the
state $x$.

Consider two cut paths from $P$ to adjacent endpoints $Q_k$
and $Q_{k+1}$.  If we move the cut paths as close together as possible
using H-moves, then by a curvature argument~\cite{Kuperberg:spiders},
there are only three possibilities for the portion of the web
between them:
$$
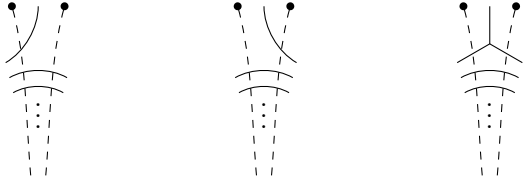
\pspicture(-1,-2)(1,1)
\pccurve[arrows=*-,linestyle=dashed,angleA=285,%
    angleB=95](-.35,.5)(-.1,-1.75)
\pccurve[arrows=*-,linestyle=dashed,angleA=255,%
    angleB=85](.35,.5)(.1,-1.75)
\pccurve[angleA=270,angleB=30](.5;90)(.5;210)
\pcarc[arcangle=30](-.383,-.45)(.383,-.45)
\pcarc[arcangle=30](-.333,-.65)(.333,-.65)
\rput(0,-.85){$\vdots$}
\endpspicture\hspace{1cm}
\pspicture(-1,-2)(1,1)
\pccurve[arrows=*-,linestyle=dashed,angleA=285,%
    angleB=95](-.35,.5)(-.1,-1.75)
\pccurve[arrows=*-,linestyle=dashed,angleA=255,%
    angleB=85](.35,.5)(.1,-1.75)
\pccurve[angleA=270,angleB=150](.5;90)(.5;330)
\pcarc[arcangle=30](-.383,-.45)(.383,-.45)
\pcarc[arcangle=30](-.333,-.65)(.333,-.65)
\rput(0,-.85){$\vdots$}
\endpspicture\hspace{1cm}
\pspicture(-1,-2)(1,1)
\pccurve[arrows=*-,linestyle=dashed,angleA=285,%
    angleB=95](-.35,.5)(-.1,-1.75)
\pccurve[arrows=*-,linestyle=dashed,angleA=255,%
    angleB=85](.35,.5)(.1,-1.75)
\psline(.5;90)(0,0)(.5;210)
\psline(.5;330)(0,0)
\pcarc[arcangle=30](-.383,-.45)(.383,-.45)
\pcarc[arcangle=30](-.333,-.65)(.333,-.65)
\rput(0,-.85){$\vdots$}
\endpspicture
$$
In each case, there are no edges with flow between the
two cut paths.  At the same time, when one performs an H-move
on a cut path, the path hops over a single edge which does not carry flow.
Thus, the set of all cut paths intersects every edge with flow.

It remains to show that the weight of the state $x$ is 1.  This
follows from comparing the extended growth rules in Figure~(\ref{fextended})
to the weights in Figures~(\ref{ftweights}) and (\ref{fbweights}).
\end{proof}

Given a sign string $S$ and a non-dominant state string $J$, the 
growth rules produce a new sign string $S'$, a state string $J'$,
and a web $w^S_J \in \Hom(V^{S'},V^S)$.  Theorem~\ref{thleading}
generalizes to the vectors $w^S_J(e^{S'}_{J'})$ to produce
a web basis for all of $V^S$ \cite{Khovanov}.

\section{The dual canonical axioms}

In this section, we will give axioms that uniquely determine the dual
canonical bases of invariant spaces.  The first axiom
involves
a certain operator $\Thetabar \in U_q(\sl(3)) \hat{\tensor}
U_q(\sl(3))$, where ``$\hat{\tensor}$'' is a certain topological
tensor product \cite{FK:canonical,Lusztig:book}.  For each tensor
product $V = V^S$, we define a $v$-antilinear endomorphism
$\Phi = \Phi^S$ inductively by the rule
$$\Phi^{SS'}(e^S \tensor e^{S'}) = \Thetabar(\Phi^S(e^S) \tensor
\Phi^{S'}(e^{S'})),$$
where $e^S \in V^S$ and $e^{S'} \in V^{S'}$.  (By $v$-antilinearity, we mean
that $\Phi$ is $\C$-linear and that $\Phi(ve) = v^{-1}\Phi(e)$.) The action
of $\Phi$ on $V^+$ and $V^-$ is defined by the stipulation that it fixes
$\{v^{\pm}_i\}$.  Remarkably, the properties of $\Thetabar$ imply that this
definition is consistent.

\begin{theorem}[Lusztig] For any sign string $S = s_1,\ldots,s_n$ and
any state string $J = j_1,\ldots,j_n$, there is a unique element
$$e^S_{\heart J} = e^{s_1}_{j_1} \heart 
e^{s_2}_{j_2} \heart \ldots \heart e^{s_n}_{j_n}
\in V^S$$
which is invariant under $\Phi$ and such that
$$e^{s_1}_{j_1} \heart e^{s_2}_{j_2} \heart \ldots \heart e^{s_n}_{j_n}
= e^S_J + \sum_{J'} c(S,J,J') e^S_{J'}$$
with $c(S,J,J') \in v^{-1} \Z[v^{-1}]$ (the negative-exponent property).
\label{thlusztig}
\end{theorem}

Clearly, $\{e^S_{\heart J}\}$ is a basis of $V^S$, the dual canonical
basis.  It is less clear, but nevertheless true, that the subset of
$\{e^S_{\heart J}\}$ indexed by dominant paths is a basis of
$\Inv(V^S)$ \cite[Sec. 27.2.5]{Lusztig:book}.

It remains to determine when $w^S_J = e^S_{\heart J}$, \ie, when basis
webs are dual canonical.

\section{Early agreement}

By inspection, the empty web $w^{\emptyset}_{\emptyset} = 1$ and the
webs $w^{-+}_{1,-1} = b^{-+}$ and $w^{+-}_{1,-1} = b^{+-}$ have the
negative-exponent property.  They are therefore dual canonical, because
there must be a dual canonical vector in the one-dimensional space of
invariants in which they lie, and they do not retain the
negative-exponent property after rescaling.  The same argument applies
to $t^{+++}$ and $t^{---}$.

\begin{proposition} Every basis web $w^S_J$ is invariant under $\Phi$.
\end{proposition}
\begin{proof}
We will actually prove that every morphism made out of $b$'s, $t$'s, and
$\sigma$'s, in other words every (monomial) web interpreted arbitrarily
as an element of $\Hom(V^S,V^{S'})$, is equivariant under $\Phi$.  For
this purpose, it is convenient to define $\Phi = \Phi^{\emptyset}$ for
a 0-fold tensor product as $v$-conjugation; $\Phi(v) = v^{-1}$.

Clearly, the identity $I$ is equivariant.

Since the $b$'s and the $t$'s are dual canonical, they are invariant
under $\Phi$, or equivariant as morphisms.  Let us assume for a moment
that the $\sigma$'s are equivariant also.  If $L$ and $L'$ are both
equivariant under both $\Phi$, then so is their composition if they can
be composed.  If they are also both equivariant under the action of
$U_q(\sl(3))$, as any web is, then $L \tensor L'$ intertwines $\Phi
\tensor \Phi$, and it also intertwines $\Thetabar \in U_q(\sl(3))
\hat{\tensor} U_q(\sl(3))$. It therefore intertwines $\Phi$.  The
proposition follows by induction, decomposing an arbitrary web as a
tensor product or composition of simpler pieces.

The equivariance of $\sigma$'s follows from equations~(\ref{eisotopy})
and a reversal of the previous argument.  The map $\Phi$ intertwines
the identity; pushing $\Phi$ from right to left on the left side of
the first equation, we conclude that
$$(I \tensor (\sigma_{-+} \circ \Phi^{-+})) \circ (b^{+-} \tensor I)
= (I \tensor (\Phi^{\emptyset} \circ \sigma_{-+})) \circ (b^{+-} \tensor I).$$
Because $b^{+-}$ is non-singular, this implies
$$\sigma_{-+} \circ \Phi^{-+} = \Phi^{\emptyset} \circ \sigma_{-+}.$$
The same argument applies to $\sigma_{+-}$.
\end{proof}

For a general basis web $w^S_J$, each state has some weight $v^n$; call
$n$ the exponent of the state.  The web $w^S_J$ has a distinguished
state with weight 1, namely the unique state with boundary $J$.  By 
Theorem~\ref{thleading}, we can call this state the leading state.  It is
dual canonical if and only if every non-leading state has negative exponent.
In the following discussion, let $w$ be a basis web which is not dual
canonical and which has as few endpoints as possible.

\begin{proposition} A minimal counterexample $w$ is connected
and does not have a Y or a double H at the boundary:
$$
w \ne 
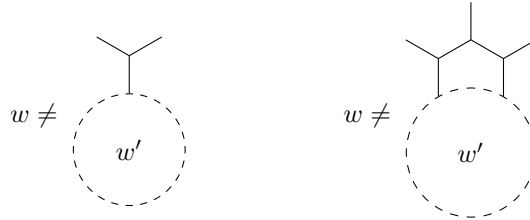
\pspicture[.47](-.85,-2.1)(.85,.5)
\psline(0,0)(.5;270)\psline(0,0)(.5;30)\psline(0,0)(.5;150)
\pscircle[linestyle=dashed](0,-1.25){.75}
\rput(0,-1.25){$w'$}
\endpspicture \hspace{2cm}
w \ne \pspicture[.47](-.966,-2.5)(.966,.5)
\psline(-.866,0)(-.433,-.25)(0,0)(.433,-.25)(.866,0)
\psline(0,0)(0,.5)\psline(-.433,-.25)(-.433,-.75)
\psline(.433,-.25)(.433,-.75)
\pscircle[linestyle=dashed](0,-1.5){.866}
\rput(0,-1.5){$w'$}
\endpspicture
$$
\label{pyh}
\end{proposition}
\begin{proof} Suppose, to the contrary, that $w$ is the disjoint union of $w'$
and $w''$.  A state of $w$ restricts to a state of $w'$ and a state of $w''$,
and the weights multiply.  Since only one state of either $w'$ or $w''$ has
weight $v^n$ with $n \ge 0$, the same is true of $w$. Thus, $w$ has
the negative-exponent property.

Suppose that $w$ has a Y at the boundary.  Let $w'$ be $w$ without the Y.   By
inspection of $t^{++}_-$ and $t^{--}_+$, the Y has six possible states, three
with weight 1 and three with weight $v^{-1}$. At the same time, any state of
$w'$ has negative exponent except for the leading state, which has exponent
$0$.  Thus, the only way that $w$ can fail to have the negative-exponent
property is if there are two distinct states which have weight 1 on both $w'$
and the Y.  But since the three states of the Y have different states on the
edge that it shares with $w'$, this is not possible.

Suppose finally that $w$ has a double H.  The argument here is essentially the
same as for a Y, but more complicated.  We arrange the double H as a
composition of one $\lambda$ and two Y's, as above.  Then one can check on a
case-by-case basis that there it has no states with positive exponent and
twelve states with weight 1. These twelve states run through all nine
possibilities for the states of the two bottom edges, with three repeats. The
three repeats are $1$ and $0$; $1$ and $-1$; and $0$ and $-1$. Again let $w'$
be the remainder. If $w$ had two states with weight $1$, then because $w'$
has the negative-exponent property, the two states can only differ in the
double H. Moreover, the states at the bottom extend to the leading state of
$w'$, so $w'$ can be constructed using the growth rules.  In particular, we
can apply a growth rule to the two edges of the double H to conclude that
$w'$ has a Y here.  Together, they make a square:
$$
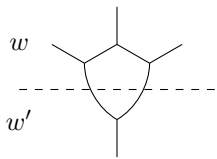
\pspicture(-1.3,-1.7)(1.3,.5)
\psline(-.866,0)(-.433,-.25)(0,0)(.433,-.25)(.866,0)
\psline(0,0)(0,.5)
\pccurve[angleA=270,angleB=150](-.433,-.25)(0,-1)
\pccurve[angleA=270,angleB= 30]( .433,-.25)(0,-1)
\qline(0,-1)(0,-1.5)
\psline[linestyle=dashed](-1.3,-.6)(1.3,-.6)
\rput(-1.3,0){$w$}
\rput(-1.3,-1){$w'$}
\endpspicture
$$
The square is an elliptic face of $w$, a contradiction.
\end{proof}

Proposition~\ref{pyh} establishes the positive part of
Theorem~\ref{thmain}, since by a curvature argument, a connected basis
web with at least one vertex and with no Y's at the boundary must have
at least six H's, and no two of the H's can share vertices.  Moreover,
one can eliminate all but one web with twelve vertices, a web which as
it happens is the first counterexample.

\section{Counterexamples}

It is easier to demonstrate that counterexamples eventually appear
than to find or verify any particularly small counterexamples.
Consider a hexagon made of three $\lambda$'s and three Y's
that is part of some larger web:
$$
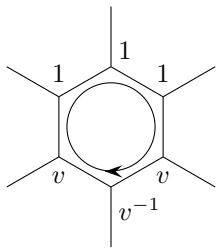
\pspicture(-1.8,-1.8)(1.8,1.8)
\pnode(.8; 90){a1}\pnode(.8;150){a2}\pnode(.8;210){a3}
\pnode(.8;270){a4}\pnode(.8;330){a5}\pnode(.8; 30){a6}
\pnode(1.6; 90){b1}\pnode(1.6;150){b2}\pnode(1.6;210){b3}
\pnode(1.6;270){b4}\pnode(1.6;330){b5}\pnode(1.6; 30){b6}
\ncline{a1}{a2}\ncline{a2}{a3}\ncline{a3}{a4}
\ncline{a4}{a5}\ncline{a5}{a6}\ncline{a6}{a1}
\ncline{a1}{b1}\ncline{a2}{b2}\ncline{a3}{b3}
\ncline{a4}{b4}\ncline{a5}{b5}\ncline{a6}{b6}
\pscircle(0,0){.593}
\psline[arrows=->,arrowscale=1.5]
    (-2.2pt,-.593)(-2.3pt,-.593)
\rput[tl](.1,-.95){$v^{-1}$}\rput[bl](.1,.9){$1$}
\rput[t](-.693,-.6){$v$}\rput[t](.693,-.6){$v$}
\rput[b](.693,.6){$1$}\rput[b](-.693,.6){$1$}
\endpspicture
$$
If the hexagon has a closed, clockwise flow line, as indicated,
then its overall weight is $v$.

A large flat web locally resembles the tiling of the plane by regular
hexagons.  The number of such hexagons can grow quadratically in the
length of the periphery, and we can put flow loops on one third of the
hexagons to form a valid state.  In the limit, the exponent of any
such state must be positive.  Thus, not all basis webs have the
negative exponent property.

The smallest counterexample is similar to the asymptotic ones:
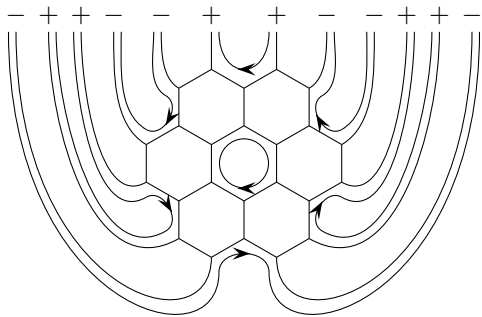
\begin{equation}
\pspicture(-3,-2.4)(3,2)
\psline(-1.299,.25)(-.866,.5)(-.433,.25)%
    (0,.5)(.433,.25)(.866,.5)(1.299,.25)
\psline(-1.299,-.25)(-.866,-.5)(-.433,-.25)(0,-.5)(.433,-.25)%
    (.866,-.5)(1.299,-.25)
\psline(-.866,1)(-.433,1.25)(0,1)(.433,1.25)(.866,1)
\psline(-.866,-1)(-.433,-1.25)(0,-1)(.433,-1.25)(.866,-1)
\psline(-.433,1.75)(-.433,1.25)\psline(.433,1.75)(.433,1.25)
\psline(-.866,1)(-.866,.5)\psline(0,1)(0,.5)\psline(.866,1)(.866,.5)
\psline(-1.299,.25)(-1.299,-.25)\psline(-.433,.25)(-.433,-.25)
\psline(1.299,.25)(1.299,-.25)\psline(.433,.25)(.433,-.25)
\psline(-.866,-1)(-.866,-.5)\psline(0,-1)(0,-.5)\psline(.866,-1)(.866,-.5)
\psbezier(-3.031,1.75)(-3.031,-1.75)(-.433,-2.5)(-.433,-1.25)
\pccurve[angleA=270,angleB=210,ncurv=.8](-2.598,1.75)(-.866,-1)
\pccurve[angleA=270,angleB=210](-2.165,1.75)(-1.299,-.25)
\pccurve[angleA=270,angleB=150](-1.732,1.75)(-1.299,.25)
\pccurve[angleA=270,angleB=150](-1.1,1.75)(-.866,1)
\pccurve[angleA=270,angleB= 30](1.1,1.75)(.866,1)
\pccurve[angleA=270,angleB= 30](1.732,1.75)(1.299,.25)
\pccurve[angleA=270,angleB=330](2.165,1.75)(1.299,-.25)
\pccurve[angleA=270,angleB=330,ncurv=.8](2.598,1.75)(.866,-1)
\psbezier(3.031,1.75)(3.031,-1.75)(.433,-2.5)(.433,-1.25)
\rput[b](-3.031,1.85){$-$}\rput[b](-2.598,1.85){$+$}
\rput[b](-2.165,1.85){$+$}\rput[b](-1.732,1.85){$-$}
\rput[b](-1.1,1.85){$-$}\rput[b](-.433,1.85){$+$}
\rput[b](.433,1.85){$+$}\rput[b](1.1,1.85){$-$}
\rput[b](1.732,1.85){$-$}\rput[b](2.165,1.85){$+$}
\rput[b](2.598,1.85){$+$}\rput[b](3.031,1.85){$-$}
\pccurve[angleA=270,angleB=330,ncurv=.7](2.498,1.75)(1.039,-.9)
\pccurve[angle=150,ncurv=1](1.039,-.9)(1.299,-.45)
\redmiddlearrow
\pccurve[angleA=270,angleB=330](2.265,1.75)(1.299,-.45)
\pccurve[angleA=270,angleB= 30](1.632,1.75)(1.299,.45)
\pccurve[angle=210,ncurv=1](1.039,.9)(1.299,.45)
\unredmiddlearrow
\pccurve[angleA=270,angleB= 30](1.2,1.75)(1.039,.9)
\pccurve[angleA=270,angleB=210,ncurv=.7](-2.498,1.75)(-1.039,-.9)
\pccurve[angle= 30,ncurv=1](-1.039,-.9)(-1.299,-.45)
\unredmiddlearrow
\pccurve[angleA=270,angleB=210](-2.265,1.75)(-1.299,-.45)
\pccurve[angleA=270,angleB=150](-1.632,1.75)(-1.299,.45)
\pccurve[angle=330,ncurv=1](-1.039,.9)(-1.299,.45)
\redmiddlearrow
\pccurve[angleA=270,angleB=150](-1.2,1.75)(-1.039,.9)
\psbezier(3.131,1.75)(3.131,-1.95)(.333,-2.7)(.333,-1.45)
\pccurve[angle=90,ncurv=1](-.333,-1.45)(.333,-1.45)\redmiddlearrow
\psbezier(-3.131,1.75)(-3.131,-1.95)(-.333,-2.7)(-.333,-1.45)
\pscircle(0,0){.333}
\psline[arrows=->,arrowscale=1.5]%
    (-2.2pt,-.333)(-2.3pt,-.333)
\pccurve[ncurv=2,angle=270](-.333,1.75)(.333,1.75)
\unredmiddlearrow
\endpspicture \label{fsmallest}
\end{equation}
One can compute the weight of this state by noting that, besides the
complete hexagon in the middle, the flow lines make three more hexagons
divided into halves, and in addition there are two right-ward pointing
U-turns.  Thus the weight is $v^4v^{-4} = 1$.  If one rotates the
web, another state with weight 1 appears:
$$
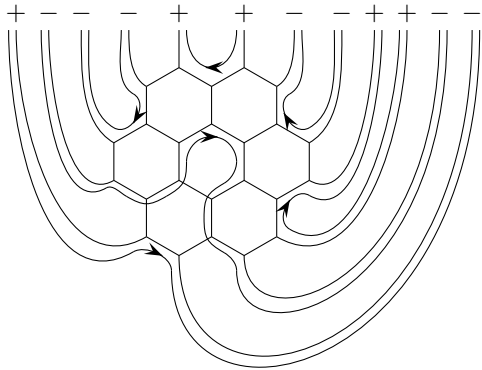
\pspicture(-2.8,-3)(3.5,2)
\psline(-1.299,.25)(-.866,.5)(-.433,.25)
    (0,.5)(.433,.25)(.866,.5)(1.299,.25)
\psline(-1.299,-.25)(-.866,-.5)(-.433,-.25)(0,-.5)(.433,-.25)
    (.866,-.5)(1.299,-.25)
\psline(-.866,1)(-.433,1.25)(0,1)(.433,1.25)(.866,1)
\psline(-.866,-1)(-.433,-1.25)(0,-1)(.433,-1.25)(.866,-1)
\psline(-.433,1.75)(-.433,1.25)\psline(.433,1.75)(.433,1.25)
\psline(-.866,1)(-.866,.5)\psline(0,1)(0,.5)\psline(.866,1)(.866,.5)
\psline(-1.299,.25)(-1.299,-.25)\psline(-.433,.25)(-.433,-.25)
\psline(1.299,.25)(1.299,-.25)\psline(.433,.25)(.433,-.25)
\psline(-.866,-1)(-.866,-.5)\psline(0,-1)(0,-.5)\psline(.866,-1)(.866,-.5)
\psbezier(3.464,1.75)(3.464,-2.75)(-.433,-3.75)(-.433,-1.25)
\pccurve[angleA=270,angleB=210,ncurv=.8](-2.598,1.75)(-.866,-1)
\pccurve[angleA=270,angleB=210](-2.165,1.75)(-1.299,-.25)
\pccurve[angleA=270,angleB=150](-1.732,1.75)(-1.299,.25)
\pccurve[angleA=270,angleB=150](-1.1,1.75)(-.866,1)
\pccurve[angleA=270,angleB= 30](1.1,1.75)(.866,1)
\pccurve[angleA=270,angleB= 30](1.732,1.75)(1.299,.25)
\pccurve[angleA=270,angleB=330](2.165,1.75)(1.299,-.25)
\pccurve[angleA=270,angleB=330,ncurv=.8](2.598,1.75)(.866,-1)
\psbezier(3.031,1.75)(3.031,-1.75)(.433,-2.5)(.433,-1.25)
\pnode(-2.698,1.75){a1}\pnode(-1.039,-1.2){a2}\pnode(-.533,-1.45){a3}
\pnode(-.533,-3.85){a4}\pnode(3.564,-2.95){a5}\pnode(3.564,1.75){a6}
\pccurve[angleA=270,angleB=210,ncurv=.8](a1)(a2)
\pccurve[angleA=30,angleB=90](a2)(a3)
\redmiddlearrow
\psbezier(a3)(a4)(a5)(a6)
\pnode(-1.083,-.375){b0}\pnode(-.65,-.375){b1}\pnode(-.433,0){b2}
\pnode(-.217,.375){b3}\pnode(.217,.375){b4}\pnode(.433,0){b5}
\pnode(.217,-.375){b6}\pnode(0,-.75){b7}\pnode(.217,-1.125){b8}

\pccurve[angleA=270,angleB=210](-2.265,1.75)(-1.399,-.35)
\pccurve[angleA= 30,angleB=150](-1.399,-.35)%
   ([angle=240,nodesep=.1]b0)
\pccurve[angleA=330,angleB=210,ncurv=1]%
   ([angle=240,nodesep=.1]b0)([angle=300,nodesep=.1]b1)
\pccurve[angleA= 30,angleB=270,ncurv=1]%
   ([angle=300,nodesep=.1]b1)([angle=  0,nodesep=.1]b2)
\pccurve[angleA= 90,angleB=210]([angle= 0,nodesep=.1]b2)%
   ([angle=300,nodesep=.1]b3)
\pccurve[angleA= 30,angleB=150]([angle=300,nodesep=.1]b3)%
   ([angle=240,nodesep=.1]b4)\redmiddlearrow
\pccurve[angleA=330,angleB= 90]([angle=240,nodesep=.1]b4)%
   ([angle=180,nodesep=.1]b5)
\pccurve[angleA=270,angleB= 30]([angle=180,nodesep=.1]b5)%
   ([angle=120,nodesep=.1]b6)
\pccurve[angleA=210,angleB= 90,ncurv=1]%
   ([angle=120,nodesep=.1]b6)([angle=180,nodesep=.1]b7)
\pccurve[angleA=270,angleB=150,ncurv=1]%
   ([angle=180,nodesep=.1]b7)([angle=240,nodesep=.1]b8)
\pccurve[angleA=330,angleB=90]([angle=240,nodesep=.1]b8)%
   (.333,-1.45)
\psbezier(3.131,1.75)(3.131,-1.95)(.333,-2.7)(.333,-1.45)

\rput[b](-2.598,1.85){$+$}
\rput[b](-2.165,1.85){$-$}\rput[b](-1.732,1.85){$-$}
\rput[b](-1.1,1.85){$-$}\rput[b](-.433,1.85){$+$}
\rput[b](.433,1.85){$+$}\rput[b](1.1,1.85){$-$}
\rput[b](1.732,1.85){$-$}\rput[b](2.165,1.85){$+$}
\rput[b](2.598,1.85){$+$}\rput[b](3.031,1.85){$-$}
\rput[b](3.464,1.85){$-$}

\pccurve[angleA=270,angleB=330,ncurv=.7](2.498,1.75)(1.039,-.9)
\pccurve[angle=150,ncurv=1](1.039,-.9)(1.299,-.45)
\redmiddlearrow
\pccurve[angleA=270,angleB=330](2.265,1.75)(1.299,-.45)
\pccurve[angleA=270,angleB= 30](1.632,1.75)(1.299,.45)
\pccurve[angle=210,ncurv=1](1.039,.9)(1.299,.45)
\unredmiddlearrow
\pccurve[angleA=270,angleB= 30](1.2,1.75)(1.039,.9)

\pccurve[angleA=270,angleB=150](-1.632,1.75)(-1.299,.45)
\pccurve[angle=330,ncurv=1](-1.039,.9)(-1.299,.45)
\redmiddlearrow
\pccurve[angleA=270,angleB=150](-1.2,1.75)(-1.039,.9)
\pccurve[ncurv=2,angle=270](-.333,1.75)(.333,1.75)
\unredmiddlearrow
\endpspicture
$$
The fact that neither web is dual canonical is consistent
with the fact that both web bases and dual canonical bases
are preserved by cyclic permutation of tensor factors.

\begin{theorem} The web
$$
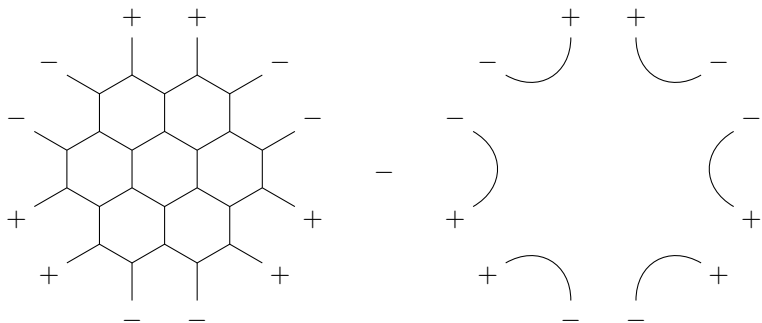
\pspicture[.47](-2.2,-2.2)(2.2,2.2)
\psline(-1.299,-1.25)(-.866,-1)(-.433,-1.25)
    (0,-1)(.433,-1.25)(.866,-1)(1.299,-1.25)
\psline(-1.732,-.5)(-1.299,-.25)(-.866,-.5)(-.433,-.25)(0,-.5)(.433,-.25)
    (.866,-.5)(1.299,-.25)(1.732,-.5)
\psline(-1.732,.5)(-1.299,.25)(-.866,.5)(-.433,.25)
    (0,.5)(.433,.25)(.866,.5)(1.299,.25)(1.732,.5)
\psline(-1.299,1.25)(-.866,1)(-.433,1.25)(0,1)(.433,1.25)(.866,1)(1.299,1.25)
\psline(-.433,1.75)(-.433,1.25)\psline(.433,1.75)(.433,1.25)
\psline(-.866,1)(-.866,.5)\psline(0,1)(0,.5)\psline(.866,1)(.866,.5)
\psline(-1.299,.25)(-1.299,-.25)\psline(-.433,.25)(-.433,-.25)
\psline(1.299,.25)(1.299,-.25)\psline(.433,.25)(.433,-.25)
\psline(-.866,-1)(-.866,-.5)\psline(0,-1)(0,-.5)\psline(.866,-1)(.866,-.5)
\psline(-.433,-1.75)(-.433,-1.25)\psline(.433,-1.75)(.433,-1.25)
\rput[b](-.433,1.9){$+$}\rput[b](.433,1.9){$+$}
\rput[br](-1.399,1.3){$-$}\rput[bl](1.399,1.3){$-$}
\rput[br](-1.832,.55){$-$}\rput[bl](1.832,.55){$-$}
\rput[tr](-1.832,-.55){$+$}\rput[tl](1.832,-.55){$+$}
\rput[tr](-1.399,-1.3){$+$}\rput[tl](1.399,-1.3){$+$}
\rput[t](-.433,-1.9){$-$}\rput[t](.433,-1.9){$-$}
\endpspicture
\hspace{.5cm} - \hspace{.5cm}
\pspicture[.47](-2.2,-2.2)(2.2,2.2)
\pccurve[angleA=270,angleB=330,ncurv=1](-.433,1.75)(-1.299,1.25)
\pccurve[angleA=330,angleB= 30,ncurv=1](-1.732,.5)(-1.732,-.5)
\pccurve[angleA= 30,angleB= 90,ncurv=1](-1.299,-1.25)(-.433,-1.75)
\pccurve[angleA= 90,angleB=150,ncurv=1](.433,-1.75)(1.299,-1.25)
\pccurve[angleA=150,angleB=210,ncurv=1](1.732,-.5)(1.732,.5)
\pccurve[angleA=210,angleB=270,ncurv=1](1.299,1.25)(.433,1.75)
\rput[b](-.433,1.9){$+$}\rput[b](.433,1.9){$+$}
\rput[br](-1.399,1.3){$-$}\rput[bl](1.399,1.3){$-$}
\rput[br](-1.832,.55){$-$}\rput[bl](1.832,.55){$-$}
\rput[tr](-1.832,-.55){$+$}\rput[tl](1.832,-.55){$+$}
\rput[tr](-1.399,-1.3){$+$}\rput[tl](1.399,-1.3){$+$}
\rput[t](-.433,-1.9){$-$}\rput[t](.433,-1.9){$-$}
\endpspicture
$$
is dual canonical. \label{thcorrection}
\end{theorem}

\begin{proof} (Sketch) Removing an H from the counterexample $w$ in
Figure~(\ref{fsmallest}) produces a non-counterexample $w'$ by
Proposition~\ref{pyh}. A hypothetical state $x$ with non-negative exponent
must either have weight $v$ or $1$ on the H; in the latter case, the state
$x$ must restrict to the leading state of $w'$ and yet differ from the
leading state of $w$.  A combinatorial investigation reveals that the second
alternative is impossible. The only possibility for the local state for each
of the six H's forces $x$ to be the state given in Figure~(\ref{fsmallest}).
Thus we can subtract off another basis web, which happens to be dual
canonical, to eliminate this term and recover the negative-exponent property.
\end{proof}

\section{What is to be done?}

What can one conclude from the fact that the web bases are not dual
canonical?  First, we argue that the web bases are nevertheless
interesting.  They are useful for computing quantum link invariants,
and they may be useful for computing $6j$-symbols along the lines of
Masbaum and Vogel's computation for $\sl(2)$ \cite{Masbaum-Vogel}.

In light of Theorem~\ref{thcorrection}, it is possible that
web bases can somehow be understood using perverse sheaves
that are the same as those related to canonical bases except
for certain salient subsheaves or stalks.

\section{Higher-rank spiders}

The main open problem concerning spiders and web bases is to generalize
the combinatorial rank two spiders to higher rank.  If the web bases
were dual canonical, it would have given an immediate definition,
albeit a very different one from the one given in
Reference~\citen{Kuperberg:spiders}.  Since it is not the case,
we outline a possible alternative approach to such a generalization.

Calculations in rank 2 spiders exhibit many elements of the
Coxeter-Weyl geometry of the corresponding affine Weyl group, this
phenomenon depends on the coincidence that the dimension of a web
equals the rank of the Lie algebra.  Nevertheless, it is implausible
that a higher-rank spider would involve higher-dimensional webs,
because quantum groups, irrespective of their rank, are fundamentally
related to low-dimensional geometry and topology.  Moreover, the $A_1$
or Temperley-Lieb spider has 2-dimensional webs even though the Lie
algebra has rank 1.

Informally, a large flat web in any of the three rank 2 spiders
resembles the Voronoi tiling of the plane associated to the
weight lattice.  More specifically, if one generically immerses a disk
in the plane of a weight lattice, then the pull-back of the
edges and vertices of the Voronoi tiling forms a valid
non-elliptic web:
$$
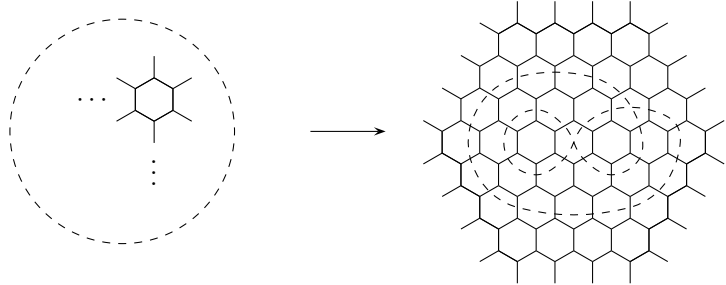
\pspicture[.5](-2,-2)(2,2)
\pscircle[linestyle=dashed](0,0){1.5}
\pnode(.6;45){a}
\rput([nodesep=.2887,angle=30]a){\littley}
\rput([nodesep=.2887,angle=90]a){\littlelam}
\rput([nodesep=.2887,angle=150]a){\littley}
\rput([nodesep=.2887,angle=210]a){\littlelam}
\rput([nodesep=.2887,angle=270]a){\littley}
\rput([nodesep=.2887,angle=330]a){\littlelam}
\rput([nodesep=.866,angle=270]a){$\vdots$}
\rput([nodesep=.8,angle=180]a){$\cdots$}
\endpspicture\psgoesto
\pspicture[.5](-2,-2)(2,2)
\rput(-1,1.299){\littley}\rput(-.5,1.299){\littley}
\rput(0,1.299){\littley}\rput(.5,1.299){\littley}
\rput(1,1.299){\littley}
\rput(-1.25,.866){\littley}\rput(-.75,.866){\littley}
\rput(-.25,.866){\littley}\rput(.25,.866){\littley}
\rput(.75,.866){\littley}\rput(1.25,.866){\littley}
\rput(-1.5,.433){\littley}\rput(-1,.433){\littley}
\rput(-.5,.433){\littley}\rput(0,.433){\littley}
\rput(.5,.433){\littley}\rput(1,.433){\littley}
\rput(1.5,.433){\littley}
\rput(-1.75,0){\littley}\rput(-1.25,0){\littley}
\rput(-.75,0){\littley}\rput(-.25,0){\littley}
\rput(.25,0){\littley}\rput(.75,0){\littley}
\rput(1.25,0){\littley}\rput(1.75,0){\littley}
\rput(-1.5,-.433){\littley}\rput(-1,-.433){\littley}
\rput(-.5,-.433){\littley}\rput(0,-.433){\littley}
\rput(.5,-.433){\littley}\rput(1,-.433){\littley}
\rput(1.5,-.433){\littley}
\rput(-1.25,-.866){\littley}\rput(-.75,-.866){\littley}
\rput(-.25,-.866){\littley}\rput(.25,-.866){\littley}
\rput(.75,-.866){\littley}\rput(1.25,-.866){\littley}
\rput(-1,-1.299){\littley}\rput(-.5,-1.299){\littley}
\rput(0,-1.299){\littley}\rput(.5,-1.299){\littley}
\rput(1,-1.299){\littley}
\rput(-.75,-1.732){\littley}\rput(-.25,-1.732){\littley}
\rput(.25,-1.732){\littley}\rput(.75,-1.732){\littley}
\rput(-.75,1.443){\littlelam}\rput(-.25,1.443){\littlelam}
\rput(.25,1.443){\littlelam}\rput(.75,1.443){\littlelam}
\rput{120}(0,-.1443){\rput(0,.1443){
\rput(-.75,1.443){\littlelam}\rput(-.25,1.443){\littlelam}
\rput(.25,1.443){\littlelam}\rput(.75,1.443){\littlelam}}}
\rput{240}(0,-.1443){\rput(0,.1443){
\rput(-.75,1.443){\littlelam}\rput(-.25,1.443){\littlelam}
\rput(.25,1.443){\littlelam}\rput(.75,1.443){\littlelam}}}
\psccurve[linestyle=dashed]%
    (.125,.0722)(-.125,-.3608)(-.5,-.5773)(-.875,-.3608)(-.875,.0722)%
    (-.5,.2887)(-.125,.0722)(.125,-.3608)(.5,-.5773)(.875,-.3608)%
    (.875,.0722)(.625,.5052)(.25,.7217)(-.75,.7217)(-1.125,.5052)%
    (-1.375,.0722)(-1.375,-.3608)(-1.125,-.7938)(-.75,-1.0103)%
    (.75,-1.0103)(1.125,-.7938)(1.375,-.3608)(1.375,.0722)(1,.2887)%
    (.5,.2887)
\endpspicture
$$
Similarly, in the rank 1 case, there is a weight-lattice Voronoi tiling
of the line by line segments.  If one submerses a disk, the inverse
image of the endpoints of these line segments is some 1-manifold, which
is then a basis web in the $A_1$ spider.

Thus, we may hypothesize that a web in a rank $n$ spider resembles the
inverse image of a weight-lattice Voronoi tiling under an immersion of
a disk in $\R^n$.  An essential ingredient, which is present in the
rank 1 and 2 cases, is that the codimension 1 faces of such a Voronoi
tiling are labelled by fundamental irreducible representations.  More
concretely, if $\mu_1$ and $\mu_2$ are two Voronoi adjacent points in
the weight lattice, there is a unique dominant weight $\lambda$ which
is conjugate to $\mu_1-\mu_2$.  The codimension 1 face separating
$\mu_1$ from $\mu_2$ pulls back under an immersion of a disk to an
edge; this edge might then be labelled by the representation
$V(\lambda)$.

For example, the weight lattice of the Lie algebra $A_3 = \sl(4)$ is
the BCC lattice.  The Voronoi region of a lattice point is a 14-side
snub octahedron; following the convention just described, the six-sided
faces are labelled with the defining representation $V$ of $\sl(4)$ and
its dual $V^*$, while the four-side faces are labelled by the
six-dimensional representation $\bigwedge^2 V$.  The incidence of the
faces suggests vertices of the form:
$$
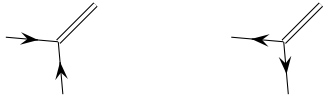
\pspicture(-1,-1)(1,1)
\pcline(.7;175)(0,0)\middlearrow\pcline(.7;275)(0,0)\middlearrow
\psline[doubleline=true](.7;45)(0,0)
\endpspicture\hspace{1cm}
\pspicture(-1,-1)(1,1)
\pcline(.7;175)(0,0)\unmiddlearrow\pcline(.7;275)(0,0)\unmiddlearrow
\pcline[doubleline=true](.7;45)(0,0)
\endpspicture
$$
Here an oriented edge is one labelled by $V$ or $V^*$,
while a double edge is one labelled by $\bigwedge^2 V$.

Both of these hypothetical vertices correspond to invariant
tensors which are unique up to a scalar factor.  One can then consider
relations which these tensors satisfy.  These include some
elliptic-looking relations such as:
$$
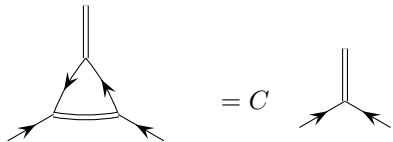
\pspicture[.45](-1.7,-1.7)(1.7,1.7)
\psline[doubleline=true](.5;90)(1.2;90)
\pcarc[arcangle=10](.5;210)(.5;90)\unmiddlearrow
\pcarc[arcangle=10](.5;90)(.5;330)\unmiddlearrow
\pcarc[arcangle=10,doubleline=true](.5;330)(.5;210)
\pcline(.5;210)(1.2;210)\unmiddlearrow
\pcline(.5;330)(1.2;330)\unmiddlearrow
\endpspicture = C \pspicture[.45](-1,-1)(1,1)
\pcline(.7;210)(0,0)\middlearrow\pcline(.7;330)(0,0)\middlearrow
\psline[doubleline=true](.7;90)(0,0)
\endpspicture
$$
But there is also the relation:
\begin{equation}
\pspicture[.45](-1.5,-1.5)(1.5,1.5)
\pcline(-.95,-.25)(-.25,-.25)\middlearrow
\pcline(-.25,-.95)(-.25,-.25)\middlearrow
\pcline(.95,.25)(.25,.25)\middlearrow
\pcline(.25,.95)(.25,.25)\middlearrow
\psline[doubleline=true](-.25,-.25)(.25,.25)
\endpspicture =
\pspicture[.45](-1.5,-1.5)(1.5,1.5)
\pcline(.95,-.25)(.25,-.25)\middlearrow
\pcline(.25,-.95)(.25,-.25)\middlearrow
\pcline(-.95,.25)(-.25,.25)\middlearrow
\pcline(-.25,.95)(-.25,.25)\middlearrow
\psline[doubleline=true](.25,-.25)(-.25,.25)
\endpspicture\label{evert}
\end{equation}
This relation (and its dual) can perhaps be motivated by a homotopy of
a disk across a vertex of the Voronoi tiling.  And there are the
relations:
\begin{align*}
\pspicture[.45](-2,-2)(2,2)
\pcline(.7;  0)(.7; 60)\unmiddlearrow\pcline[doubleline=true](.7; 60)(.7;120)
\pcline(.7;120)(.7;180)\unmiddlearrow\pcline[doubleline=true](.7;180)(.7;240)
\pcline(.7;240)(.7;300)\unmiddlearrow\pcline[doubleline=true](.7;300)(.7;  0)
\pcline(.7;  0)(1.4;  0)\unmiddlearrow\pcline(.7; 60)(1.4; 60)\middlearrow
\pcline(.7;120)(1.4;120)\unmiddlearrow\pcline(.7;180)(1.4;180)\middlearrow
\pcline(.7;240)(1.4;240)\unmiddlearrow\pcline(.7;300)(1.4;300)\middlearrow
\endpspicture &=
\pspicture[.45](-2,-2)(2,2)
\pcline[doubleline=true](.7;  0)(.7; 60)\pcline(.7; 60)(.7;120)\middlearrow
\pcline[doubleline=true](.7;120)(.7;180)\pcline(.7;180)(.7;240)\middlearrow
\pcline[doubleline=true](.7;240)(.7;300)\pcline(.7;300)(.7;  0)\middlearrow
\pcline(.7;  0)(1.4;  0)\unmiddlearrow\pcline(.7; 60)(1.4; 60)\middlearrow
\pcline(.7;120)(1.4;120)\unmiddlearrow\pcline(.7;180)(1.4;180)\middlearrow
\pcline(.7;240)(1.4;240)\unmiddlearrow\pcline(.7;300)(1.4;300)\middlearrow
\endpspicture + \mbox{simpler terms} \\
\pspicture[.45](-1.5,-1.5)(1.5,1.5)
\pcline[doubleline=true](.5; 45)(1.2; 45)
\pcline[doubleline=true](.5;135)(1.2;135)
\pcline[doubleline=true](.5;225)(1.2;225)
\pcline[doubleline=true](.5;315)(1.2;315)
\pcline(.5; 45)(.5;135)\middlearrow\pcline(.5;135)(.5;225)\unmiddlearrow
\pcline(.5;225)(.5;315)\middlearrow\pcline(.5;315)(.5; 45)\unmiddlearrow
\endpspicture &=
\pspicture[.45](-1.5,-1.5)(1.5,1.5)
\pcline[doubleline=true](.5; 45)(1.2; 45)
\pcline[doubleline=true](.5;135)(1.2;135)
\pcline[doubleline=true](.5;225)(1.2;225)
\pcline[doubleline=true](.5;315)(1.2;315)
\pcline(.5; 45)(.5;135)\unmiddlearrow\pcline(.5;135)(.5;225)\middlearrow
\pcline(.5;225)(.5;315)\unmiddlearrow\pcline(.5;315)(.5; 45)\middlearrow
\endpspicture + \mbox{simpler terms}
\end{align*}
These relations appear to be related to the faces of the snub
octahedron.

Despite these suggestive relations, we do not know how to put them into
a coherent whole.  Some relations that one might predict from Voronoi
geometry do not hold.  For example, if one maintains that the structure
of a vertex in the Voronoi tiling predicts relation~(\ref{evert}), then
presumably it would also predict a relation between the two webs:
$$
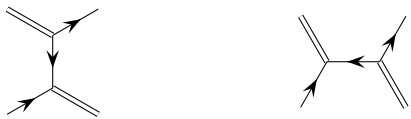
\pspicture(-1,-1)(1,1)
\pcline(0, .35)(.606, .7)\middlearrow
\pcline[doubleline=true](0, .35)(-.606, .7)
\pcline[doubleline=true](0,-.35)(.606,-.7)
\pcline(0,-.35)(-.606,-.7)\unmiddlearrow
\pcline(0,-.35)(0,.35)\unmiddlearrow
\endpspicture
\hspace{2cm}
\pspicture(-1,-1)(1,1)
\pcline( .35,0)( .7,.606)\middlearrow
\pcline[doubleline=true]( .35,0)( .7,-.606)
\pcline[doubleline=true](-.35,0)(-.7,.606)
\pcline(-.35,0)(-.7,-.606)\unmiddlearrow
\pcline(-.35,0)(.35,0)\unmiddlearrow
\endpspicture
$$
However, these two webs are linearly independent.

Ideally, we would like an explicit presentation for the representation
category of a Lie algebra which is akin to the Serre relations for the
Lie algebra itself.


\end{document}